\documentclass[11pt, a4paper]{article}
\usepackage{times}
\usepackage{a4wide}
\usepackage{hyperref}
\usepackage[british]{babel}
\usepackage{enumerate, longtable}
\usepackage{amsmath, amscd, amsfonts, amsthm, amssymb, latexsym, comment, graphicx, color, float, multicol}
\usepackage[all]{xypic}
\usepackage[T1]{fontenc}
\usepackage[utf8]{inputenc}
\usepackage{verbatim}

\usepackage{chngcntr}
\counterwithin{equation}{section}

\newtheorem{thm}{Theorem}[section]

\newtheorem{conjecture}[thm]{Conjecture}

\newtheorem{question}[thm]{Question}

\newcommand{\donothing}[1]{}

\newcommand{\Laplace}{\mathrm{Laplace}}
\newcommand{\Cauchy}{\mathrm{Cauchy}}
\newcommand{\GG}{\mathrm{GG}}

\newcommand{\exclude}[1]{}

\begin{document}

\title{On the distribution of coefficients of half-integral weight modular forms and the Bruinier-Kohnen Conjecture}
\author{Ilker Inam \footnote{I.I., E.T. Bilecik Seyh Edebali University, Department of Mathematics, Faculty of Arts and Sciences, 11200 Bilecik, Turkey, ilker.inam@bilecik.edu.tr, eliftercan5@gmail.com}, Zeynep Demirkol \"{O}zkaya \footnote{Z.D.\"{O}. Van Yuzuncu Yil University, Muradiye Vocational School, 65000 Van, Turkey, zeynepdemirkolozkaya@yyu.edu.tr} , Elif Tercan\, and Gabor Wiese\footnote{G.W. University of Luxembourg, Department of Mathematics, Maison du Nombre 6, Avenue de la Fonte, L-4364 Esch-sur-Alzette, Luxembourg, gabor.wiese@uni.lu }}
\maketitle

\begin{abstract}
This work represents a systematic computational study of the distribution of the Fourier coefficients of cuspidal Hecke eigenforms of level $\Gamma_0(4)$ and {\em half-integral} weights. Based on substantial calculations, the question is raised whether the distribution of normalised Fourier coefficients with bounded indices can be approximated by a generalised Gaussian distribution. Moreover, it is argued that the apparent symmetry around zero of the data lends strong evidence to the Bruinier-Kohnen Conjecture on the equidistribution of signs and even suggests the strengthening that signs and absolute values are distributed independently.
\end{abstract}

\section{Introduction}
This article represents a systematic computational study of the Fourier coefficients of {\em half-integral weight} cuspidal Hecke eigenforms with the aim of experimentally shedding new light on their distribution, particularly focusing on signs.
\smallskip

\noindent{\bf Size and normalisation of coefficients and the Ramanujan-Petersson Conjecture.}\\ Let $f = \sum_{n=1}^\infty a(n) q^n$ be a cuspidal Hecke eigenform of weight~$k$. The Ramanujan-Petersson Conjecture (see e.g.~\cite{30}) claims\[ a(n) = O(n^{(k-1)/2 + \epsilon})\]for any $\epsilon > 0$. By Deligne's famous proof of the Weil Conjectures~\cite{16}, the Ramanujan--Petersson Conjecture is true with $\epsilon=0$ in the {\em integer} weight case. Motivated by the Ramanujan-Petersson Conjecture we define the {\em normalised coefficients} to be\[ b(n) := \frac{a(n)}{n^{(k-1)/2}}.\] It seems that the Ramanujan-Petersson Conjecture has not been proved for a single cuspform of {\em non-integral} weight. However, it is known by work of Gun and Kohnen in \cite{18} that the Ramanujan--Petersson Conjecture would fail for $\epsilon = 0$ in half-integral weight. Their argument uses a sequence of non-squarefree indices coming from the Shimura lift to construct a counter example. In recent work of Gun, Kohnen and Soundararajan \cite{19}, the authors suggest that in half-integral weight `perhaps' the bound
\begin{equation*}
|b(|n|)| \le \exp(C \sqrt{\log|n| \log\log |n|}),
\end{equation*} derived from conjectures in~\cite{17} and stronger than the stated form of the Ramanujan-Petersson Conjecture might hold.
\smallskip

\noindent{\bf Known results and conjectures on the distribution in half-integral weight.}
In half-integral weight $k= \ell + \frac{1}{2}$, there is the crucial relation, due to Waldspurger (Theorem 1 on p.~378 of~\cite{42}) and Kohnen-Zagier (Theorem~1 of \cite{31}, p.~177), between the {\em squares} of the Fourier coefficients and central values of $L$-functions. More precisely, the Shimura lift (see \cite{38} and \cite{40}) relates those Fourier coefficients of~$f$ that are indexed by $tn^2$ with $t \in \mathbb{N}$ squarefree and $n \in \mathbb{N}$ to the $n$-th Fourier coefficient of a modular form~$g$ in integral weight~$2\ell$. Then $b(|n|)^2$ is proportional to $L(g,\chi_n,\ell)$ for fundamental discriminants $n$ such that $|n| =(-1)^\ell n$, where $L(g,\chi_n,s)$ is the Hecke L-function of~$g$ twisted by the primitive quadratic character $\chi_n$ corresponding to~$n$.

This relation is at the basis of most of the results on the absolute value of $b(|n|)$ and has led to a conjectural description of the distribution of the $b(|n|)^2$ for fundamental discriminants~$n$. In that context, we recall that the famous Sato-Tate Conjecture describing the distribution of the normalised Fourier coefficients in the {\em integral weight} case has been proved by Barnet-Lamb, Geraghty, Harris and Taylor~\cite{11}. In \cite{13}, Conrey, Keating, Rubinstein and Snaith propose a conjecture on the distribution of coefficients of modular forms of weight $3/2$ attached to elliptic curves. Their Conjecture~4.2 states that for a modular form of weight $3/2$ which is attached to an elliptic curve~$E$, the natural density of fundamental discriminants $n$ such that\[ \alpha \le (\kappa^\pm \sqrt{\log(|n|)}(b(|n|))^2)^{\frac{1}{\sqrt{\log \log |n|}}}\le \beta\] equals\[\frac{1}{\sqrt{2 \pi}}\int_\alpha^\beta \frac{1}{t} \exp(-\frac{1}{2} (\log t )^2) dt,\]where $\kappa^\pm$ is a positive constant and $0 \leq \alpha \leq \beta$. K.~Soundararajan kindly informed us that similar conjectures are made for higher weights as well. Note also that Conjecture 4.2 of~\cite{13} implies that the normalised coefficients $b(|n|)$ tend to zero almost surely, in the sense that for all $\epsilon>0$, the set of $d \in S^{\pm}$ such that $|b(|d|)| < \epsilon$ has natural density equal to~$1$. This prediction is confirmed by a theorem of Radziwi\l\l{} and Soundararajan~\cite{39}, as cited in~\cite{19} in level $\Gamma_0(4)$:\begin{quote}{\em For every $\epsilon>0$, there is a constant $C=C(\epsilon,f)$ such that for all but $o(x)$ fundamental discriminants $n$ with $x \le (-1)^\ell n \le 2x$, one has
\begin{equation}\label{eq:upper} |b(|n|)| \le C \cdot \log(|n|)^{-\frac{1}{4} + \epsilon}.
\end{equation}}
\end{quote}This theorem hence implies that the normalised coefficients $b(n)$ tend to~$0$ with probability~$1$ and hence follow a Dirac distribution at~$0$. In the other direction, in recent work of Gun, Kohnen and Soundararajan \cite{19}, the existence of large values for normalised Fourier is proved (for level $\Gamma_0(4)$):
\begin{quote}{\em For any $\epsilon>0$ and $x$ large, there are at least $x^{1-\epsilon}$ fundamental discriminants~$n$ with $x < (-1)^\ell n < 2x$ such that\begin{equation}\label{eq:lower}|b(|n|)| \ge \exp\left(\frac{1}{82}\frac{\sqrt{\log |n| }}{\sqrt{\log\log |n|}}\right).
\end{equation}}
\end{quote}

Note that the relation with central values of $L$-functions only gives information about the squares of the coefficients and hence no information about their signs. This is where the conjecture of Bruinier and Kohnen (\cite{10}, \cite{28}) enters, claiming that exactly half of the non-zero coefficients are positive, the other half negative.

At the heart of the signs problem on Fourier coefficients of modular forms of (half)-integer weights, there lie non-vanishing results such as \cite{1}, \cite{2}, \cite{3}, \cite{4} and \cite{5}. We note that there are also some interesting results on the (small) gaps between nonzero Fourier coefficients, for instance \cite{14}, \cite{15}, \cite{32} and \cite{33}. The sign change problem receives much attention and for example results \cite{12} and \cite{27} are in this direction. There are also more general results as in \cite{34}, \cite{35} and one fundamental paper in another direction is \cite{36}. Finally, there has been important recent progress on the distribution of the signs. In \cite{37}, the authors show for a Hecke eigenform in the Kohnen-plus space that $a(n)$ is negative for a positive proportion of fundamental discriminants $n$ under the assumption of the Generalised Riemann Hypothesis, as well as a similar result for positive $a(n)$. In a continuation work, in \cite{25} the authors generalize Theorem 1 of \cite{19} and they extend the results of \cite{37} to non-eigenforms.

\smallskip\noindent{\bf Contributions of this work.} Here are the main points that we want to make with this article.
\begin{enumerate}[(1)]
\item Even though it is known by \eqref{eq:upper} that the absolute values of the normalised coefficients $b(|n|)$ tend to zero with probability~$1$, we do observe a very neat non-trivial distribution of the coefficients up to (computationally reachable) bounds. The distribution seems to follow a generalised Gaussian distribution.
\item The histograms of the distribution of the normalised Fourier coefficients up to varying bounds and for varying Hecke eigenforms of half-integral weights all seem to present a similar `global shape' in the sense that they can be well approximated by a single type of density function, and that only the parameters depend on the modular form and on the bound.
\item The symmetry around $0$ of the observed distributions of the coefficients up to bounds can be interpreted as very strong evidence towards the Bruinier-Kohnen Conjecture. In fact, it suggests a strengthening of the conjecture to the point that the absolute value and the signs are independently distributed (see Conjecture~\ref{conj:SBK}). To the best of our knowledge, the calculations in this article can be seen as the most systematic and largest computational support for the Bruinier-Kohnen Conjecture to date. Furthermore, if Question \ref{qu:GGG} has a positive answer then the Bruinier-Kohnen Conjecture is true and this links the two topics in the title. 
\end{enumerate}

\smallskip\noindent{\bf Short overview over the article.} In \S\ref{sec:examples}, we report on the examples of half-integral weight modular cuspforms used for our experimental study and how they were computed. The main point is that we chose to stay in the lowest possible level $\Gamma_0(4)$ and considered weights up to $61/2$. In view of studying the distribution of the normalised Fourier coefficients of the examples, we created histograms and report on them in~\S\ref{sec:histograms}. We take into account the specific nature of half-integral weight forms that distinguishes them significantly from the well understood integral weight ones: we disregard all coefficients that via the Shimura lift come from the integral case and, consequently, study only squarefree indexed coefficients. The similar shape that the histograms exhibit suggests that they can be described by distribution functions. We consider four types of such functions in \S\ref{sec:distribution}, namely, the Laplace and the Cauchy distribution as well as two generalisations of the Gaussian distribution. We also report on data obtained when fitting the aforementioned distribution functions with the histograms. It turns out that one of the generalised Gaussian distributions is clearly the best one. In view of the fact that the normalised coefficients tend to zero almost surely, in \S\ref{sec:dependence} we explicitly seek for dependencies of the best fit parameters with data such as the number of coefficients used. We also formulate the explicit question if the normalised Fourier coefficients up to any bound indeed follow such a generalised Gaussian distribution, see Question~\ref{qu:GGG}. Finally, in \S\ref{sec:BK}, we recall the Bruinier-Kohnen Conjecture and some previous results on it. We make the point that the observed symmetry of the histograms and of the studied distribution functions around zero is strong computational evidence in favour of the conjecture and even suggests a strengthening of it.

\section{Examples of Hecke eigenforms in half-integral weights for $\Gamma_0(4)$}\label{sec:examples}
For studying their distribution computationally, we need as many Fourier coefficients of modular forms as possible. Since we are interested in higher weights, we choose to work in the smallest possible level $\Gamma_0(4)$. As described in the article \cite{24} by two of the authors, the Kohnen-plus space in half-integral weight admits bases that can be computed relatively quickly up to some high precision. For this article, we worked with the Rankin-Cohen basis as described in loc.~cit. We performed exact computations over the rationals in order not to lose any precision and only converted the normalised coefficients to real numbers in the end. A disadvantage of this choice is a huge consumption of memory, when the $q$-expansions are computed up to a high power of~$q$.

To give some more details, we use Pari/GP (see \cite{41}) to express Hecke eigenforms with respect to the Rankin-Cohen basis. Here the {\em mf package} of Pari/GP \cite{8} provides us with the necessary tools. Then, we export the basis coefficients to Magma~\cite{9}, where we construct the Hecke eigenform as a power series (in general, over a number field). This is done in Magma because it provides very fast algorithms for the multiplication of power series. In a final step, we compute the normalised coefficients over the reals. Since all previous computations are exact computations, $10$ digits of real precision are enough.

We only recorded coefficients at squarefree indices which are not known to be zero by the fact that we look only in the Kohnen-plus space (e.g.\ when $k-1/2$ is even, $a(n)$ is zero when $n$ is $2$ or $3$ modulo~$4$). We also normalised all modular forms in such a way that the first recorded normalised coefficient equals~$1$. This is the natural way of normalisation if we consider the definition of the Kohnen-plus space, but it is in no way canonical.

We compute all Hecke eigenforms of weights $13/2, 17/2, 19/2, \dots, 61/2$ (level $\Gamma_0(4)$) with $10^7$ Fourier coefficients. By this, we mean all normalised coefficients $b(n)$ with squarefree index $n < 10^7$. We reached $10^8$ Fourier coefficients for some examples and for the weight $13/2$, we could go up to $2 \cdot 10^8$. Text files containing the normalised coefficients can be downloaded\footnote{\url{http://math.uni.lu/wiese/FourierData.html}}. The total amount of data used for the study exceeds 4 GB. In all tables below, a label such as $25/2(2)$ stands for the second cuspidal Hecke eigenform (with respect to an internal ordering) in weight $25/2$ and level $\Gamma_0(4)$. The reader is referred to \cite{24} for some more details on the computations.

Note that, under the Shimura lift, any half-integral weight Hecke eigenform (in weight~$k$) corresponds to an integral weight Hecke eigenform in level~$1$ and weight $2k-1 \in 2\mathbb{Z}$. By \cite[p.~241]{29}, the Shimura lift is a Hecke equivariant isomorphism between the Kohnen plus space and the corresponding space in integral weight. This means that the eigenforms in half-integral weight are in bijection with those in integral weight. By Maeda's Conjecture (see \cite{21}), in weight $2k-1$ there is only a single Hecke orbit of eigenforms. Assuming Maeda's conjecture (which is known up to high weight by~\cite{20}, far exceeding our examples), it follows that the number of half-integral weight Hecke eigenforms in the Kohnen-plus space equals the degree of the number field generated by the coefficients of the integral weight form.

\section{Histograms for the distribution of normalised coefficients}\label{sec:histograms}
The principal aim of this article is to understand the distribution of normalised coefficients of half-integral weight Hecke eigenforms. More precisely, the point we want to make is the following. Even though the normalised coefficients $b(n)$ tend to zero almost surely by the cited result of Radziwi\l\l{} and Soundararajan, the coefficients up to bounds that we can computationally reach do follow a very interesting non trivial distribution, which can be well approximated by density functions. They turn out to be symmetric around zero.

In order to study the distribution, we created histograms for the distribution of the normalised coefficients for all the modular forms mentioned in the previous section. We restricted our attention to coefficients with squarefree indices that are not known to be zero by the fact that the modular form lies in the Kohnen plus-space. The reason for only considering squarefree indices is the following: coefficients at indices of the form $t n^2$ with $t$ squarefree and $n \in \mathbb{Z}_{\ge 2}$ are governed by the $n$-th coefficient of the Shimura lift, which is an integral weight eigenform and thus behaves with respect to the proved Sato-Tate law (if it is not CM). So, if we did not restrict to squarefree indices, we would `mix' two distributions, making the pictures harder to analyse.

\begin{figure}[h]
\centering
\includegraphics[width=\textwidth]{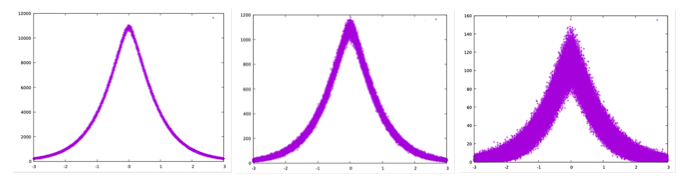}
\caption{Histogram of $10^8$ normalised Fourier coefficients of the Hecke eigenform of weight $13/2$ with box sizes $0.001$, $0.0001$ and $0.00001$, respectively}\label{figbox}
\end{figure}

We created histograms for the distribution of the normalised coefficients using gnuplot~\cite{43}. One choice that one has to make is that of the box size for the histograms. In order to understand dependencies, we created the histograms with different box sizes. Some box sizes are more pleasing to the eye than others (sometimes depending on the modular form). The graphs in Figure~\ref{figbox} are the histograms of the normalised Fourier coefficients $b(n)$ for the Hecke eigenform of weight $13/2$ with box sizes $0.001$, $0.0001$ and $0.00001$, respectively, with $10^8$ Fourier coefficients.

We observed that the choice of box size does not influence the parameters for best fits in any significant way. To confirm this, we give Table~\ref{Best} for the parameter $a$ for the $GG$-distribution for five different forms (indicated by their weights) and three different box widths. We hence disregard box sizes in our discussions and, throughout the paper, we use graphs with box size $0.001$ since in this case, the graphs seem the most pleasing to the eyes.

\begin{table}[h]\centering\begin{tabular}{|l||l|l|l|}
\hline
& 0.001 & 0.0001 & 0.00001\\
\hline
\hline
$13/2$ & 0.634 & 0.634 & 0.633\\
\hline
$17/2$ & 0.553 & 0.553 & 0.553\\
\hline
$21/2$ & 0.558 & 0.558 & 0.557\\
\hline
$25/2(1)$ & 0.504 & 0.504 & 0.504\\
\hline
$29/2(1)$ & 0.506 & 0.506 & 0.506\\
\hline
\end{tabular}\caption{Best fit parameters with different box widths for the parameter $a$ in the GG-distribution}\label{Best}
\end{table}

\section{Candidate distribution functions and fits}\label{sec:distribution}
The point that we want to make in this section is the following: The `global shape' of the histograms is independent of the modular form and of the bound for the coefficients. More precisely, our computations suggest that the histograms for the normalised Fourier coefficients of any half-integral weight Hecke eigenforms up to a given bound can be described by a single type of density function, and that only the parameters depend on the modular form as well as on the bound.

Looking at the histograms, one immediately notices that the histograms look symmetric around~$0$. Even though they present some kind of bell shape, one also sees very quickly that they do not follow a standard Gaussian. Instead, we tried the following two generalisations of the standard Gaussian and also the Laplace and the Cauchy distributions:
\begin{align}GGG(x) & := b \cdot \exp(- \frac{(d+x^2)^a }{c}) \label{eq:GGG}\\GG(x) & := b \cdot \exp(- \frac{(x^2)^a }{c}) \label{eq:GG} \\Laplace(x) & := b \cdot \exp(- |x|/c) \label{eq:Laplace} \\Cauchy(x) & := \frac{a}{b+(cx)^2}. \label{eq:Cauchy}
\end{align}

Of course, $GG$ is a special case of $GGG$ (with $d=0$) and $Laplace$ is a special case of~$GG$ (with $a=0.5$).Since the $a$-parameter in $GG$ is quite close to $0.5$ (the data is given in the appendix) and because the Laplace distribution is much simpler than the Generalised Gaussians, we took it up into our considerations. Since we did not normalise our histograms so that the area under it is~$1$, we also did not normalise the above distribution functions to be probability distributions (even though we think of them this way).

\begin{figure}[H]
\begin{center}
\includegraphics[width=0.55\textwidth]{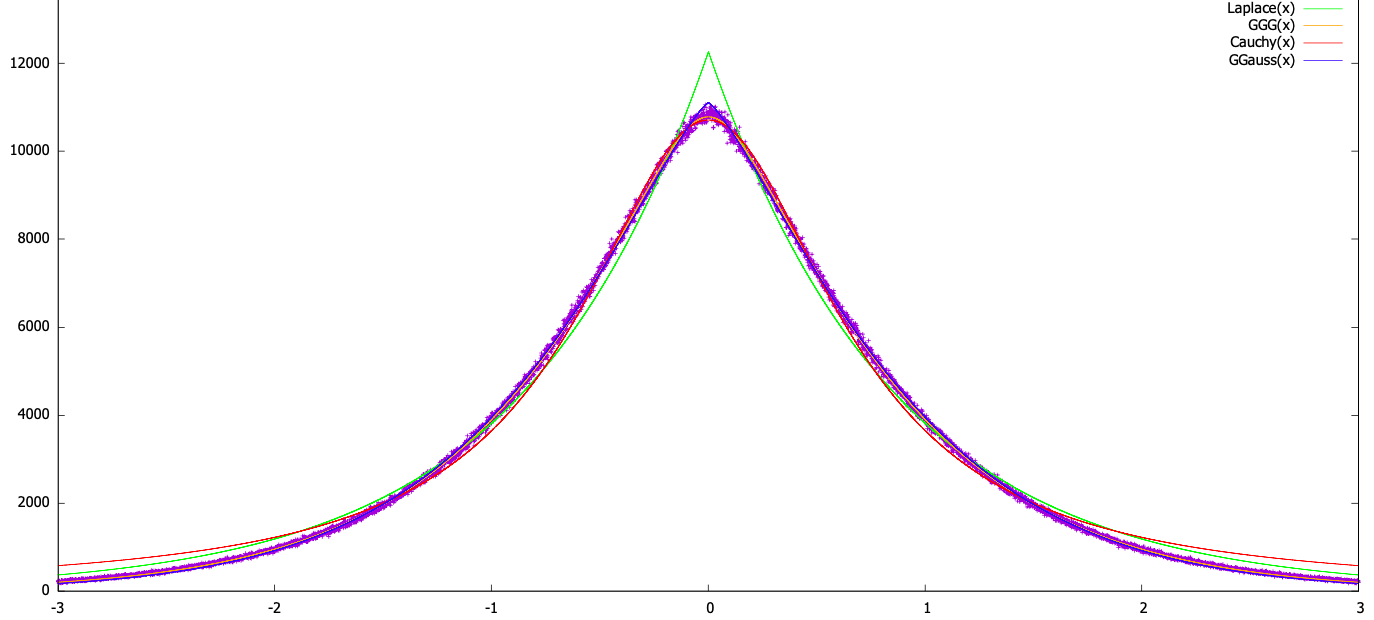}
\end{center}\caption{Histogram and distributions of Hecke eigenform of weight $13/2$ with $10^8$ coefficients}
\label{fig:1}
\end{figure}

\begin{figure}[H]
\begin{center}
\includegraphics[width=0.55\textwidth]{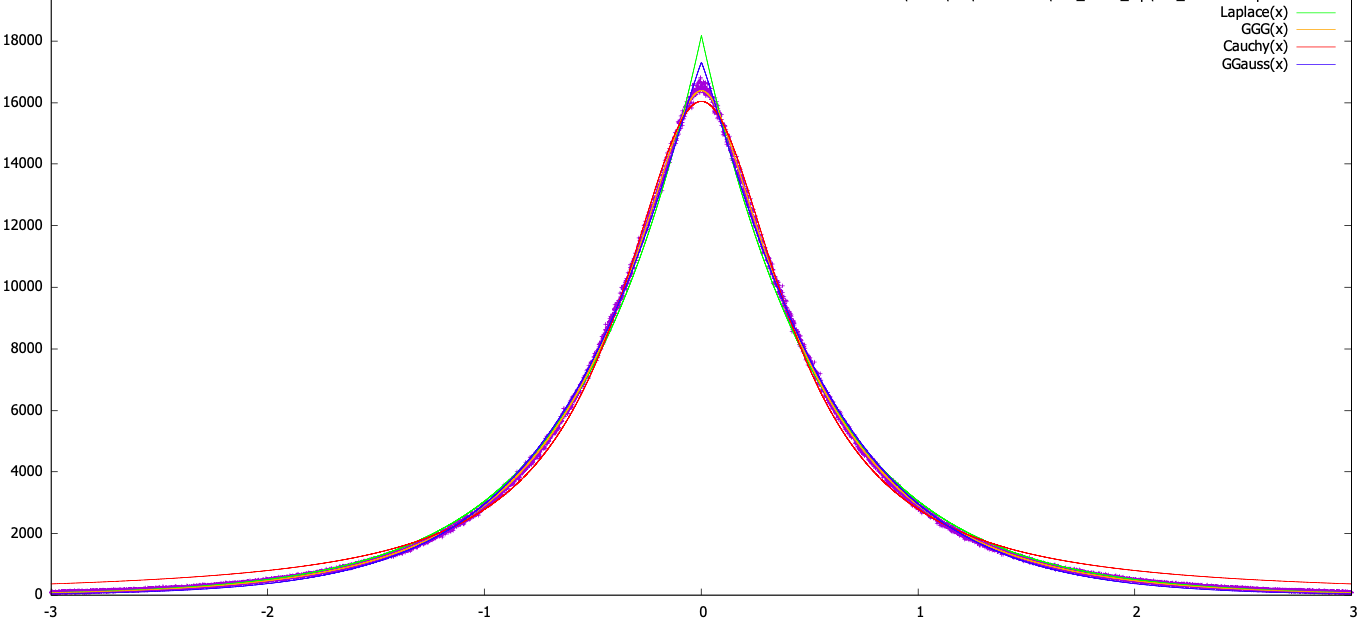}
\end{center}
\caption{Histogram and distributions of Hecke eigenform of weight $25/2$ with $10^8$ coefficients}
\label{fig:2}
\end{figure}

\begin{figure}[H]
\begin{center}
\includegraphics[width=0.55\textwidth]{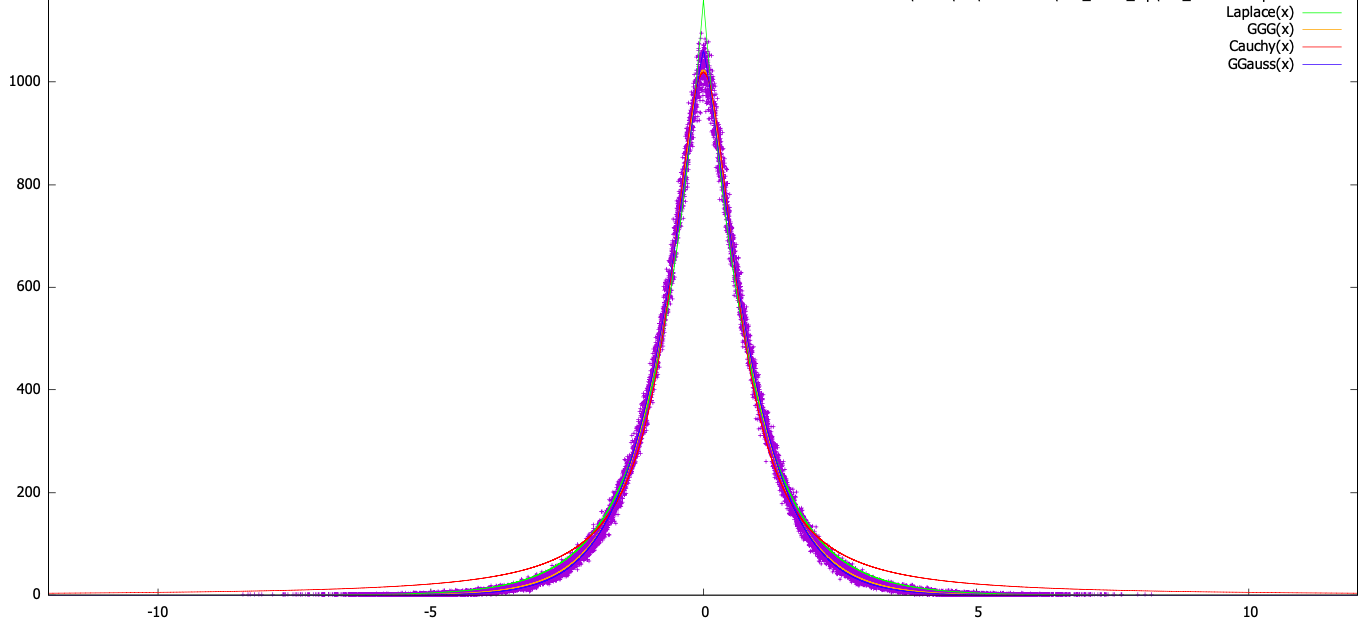}
\end{center}
\caption{Histogram and distributions of the second Hecke eigenform of weight $43/2$ with $10^7$ coefficients}
\label{fig:3}
\end{figure}

Graphically, all four distribution functions describe the histograms pretty well, the $GGG$-distribution being clearly the best. The Cauchy distribution seems to be systematically too high in the tails, whence we consider it the worst of the four. The reader is referred to the appendix for the graphs with inscribed best fit distribution functions. Sample graphs of some fits are given in Figures \ref{fig:1}, \ref{fig:2} and~\ref{fig:3}.

Since our histograms are not uniformly normalised (recall that we normalised the coefficients of the half-integral weight form in such a way that the first non-zero coefficient is~$1$) in the sense that generally they present the same shape, but someare wider, some are steeper, etc., it is very hard to compare the quality of the fits between different histograms. We will measure the quality of the fits by the Root Mean Square (RMS) value as output by gnuplot. Of course, the $GGG$-fit will always be better than the $GG$-fit and that one will always be better than the $Laplace$-fit as they are special cases of each other.

We illustrate the results of the fits performed using gnuplot by giving in Table~\ref{table:2} all best fit values for all examples for which we reached $10^8$ Fourier coefficients. The results for computations up to $10^7$ Fourier coefficients are included in the appendix.


\begin{table}
\caption{}\label{table:2}
\begin{multicols}{2}
\noindent Best fit parameters (rounded)\\ for the $GGG$-distribution:\\[.2cm]
\begin{tabular}{|l||r|r|r|r|}\hline& $a$ & $b$ & $c$ & $d$\\\hline\hline$13/2$ & 0.581 & 12538 & 0.872 & 0.030\\\hline$17/2$ & 0.453 & 20421 & 0.550 & 0.030\\\hline$19/2$ & 0.385 & 44105 & 0.317 & 0.012 \\\hline$21/2$ & 0.460 & 23411 & 0.485 & 0.022\\\hline$23/2$ & 0.494 & 14866 & 0.725 & 0.034\\\hline$25/2(1)$ & 0.391 & 19462 & 0.577 & 0.035\\\hline$25/2(2)$ & 0.237 & 88927 & 0.284 & 0.033\\\hline$27/2$ & 0.513 & 22681 & 0.471 & 0.014\\\hline$29/2(1)$ & 0.364 & 11886 & 0.812 & 0.092\\\hline$29/2(2)$ & 0.423 & 30278 & 0.402 & 0.016\\\hline
\end{tabular}

\noindent Best fit parameters (rounded)\\ for the $GG$-distribution:\\[.2cm]
\noindent 
\begin{tabular}{|l||r|r|r|}\hline& $a$ & $b$ & $c$\\\hline\hline$13/2$ &0.634 & 11105 & 0.969\\\hline$17/2$ & 0.553 & 14999 & 0.663\\\hline$19/2$ & 0.515 & 26822 & 0.363\\\hline$21/2$ & 0.558 & 17300 & 0.566\\\hline$23/2$ & 0.573 & 11997 & 0.857\\\hline$25/2(1)$ & 0.504 & 13107 & 0.752\\\hline$25/2(2)$ & 0.453 & 20676 & 0.485\\\hline$27/2$ & 0.584 & 18721 & 0.514\\\hline$29/2(1)$ & 0.506 & 7555 & 1.256\\\hline$29/2(2)$ & 0.532 & 20884 & 0.466\\\hline
\end{tabular}
\end{multicols}

\begin{multicols}{2}
\noindent Best fit parameters (rounded)\\ for the Laplace distribution:\\[.2cm]
\noindent 
\begin{tabular}{|l||r|r|}\hline& $b$ & $c$ \\\hline\hline$13/2$ & 12264 & 0.860\\\hline$17/2$ & 15711 & 0.650\\\hline$19/2$ & 27201 & 0.367\\\hline$21/2$ & 18184 & 0.560\\\hline$23/2$ & 12746 & 0.805\\\hline$25/2(1)$ & 13161 & 0.750\\\hline$25/2(2)$ & 19699 & 0.484\\\hline$27/2$ & 20071 & 0.514\\\hline$29/2(1)$ & 7597 & 1.244\\\hline$29/2(2)$ & 21515 & 0.470\\\hline
\end{tabular}

\noindent Best fit parameters (rounded)\\ for the Cauchy distribution:\\[.2cm]
\begin{tabular}{|l||r|r|r|}\hline& $a$ & $b$ & $c$\\\hline\hline$13/2$ & 183 & 0.017 & 0.181\\\hline$17/2$ & 194 & 0.014 & 0.222\\\hline$19/2$ & 2455 & 0.010 & 0.340\\\hline$21/2$ & 194 & 0.012 & 0.240\\\hline$23/2$ & 208 & 0.019 & 0.204\\\hline$25/2(1)$ & 194 & 0.017 & 0.212\\\hline$25/2(2)$ & 1376 & 0.080 & 0.725\\\hline$27/2$ & 233 & 0.013 & 0.272\\\hline$29/2(1)$ & 813 & 0.122 & 0.340\\\hline$29/2(2)$ & 342 & 0.018 & -0.350\\\hline
\end{tabular}
\end{multicols}

\noindent RMS values:
\begin{multicols}{2}\noindent 
\begin{tabular}{|l||r|r|r|r|}\hline& GG & GGG & Laplace & Cauchy \\\hline\hline$13/2$ & 76 & 59 & 286 & 252\\\hline$17/2$ & 126 & 64 & 187 & 219\\\hline$19/2$ & 208 & 75 & 215 & 295 \\\hline$21/2$ & 133 & 66 & 209 & 265\\\hline$23/2$ & 97 & 60 & 190 & 189\\\hline
\end{tabular}

\noindent 
\begin{tabular}{|l||r|r|r|r|}\hline& GG & GGG & Laplace & Cauchy \\\hline\hline$25/2(1)$ & 134 & 62 & 135 & 178\\\hline$25/2(2)$ & 283 & 79 & 320 & 238\\\hline$27/2$ & 111 & 63 & 266 & 347\\\hline$29/2(1)$ & 95 & 56 & 95 & 157\\\hline$29/2(2)$ & 161 & 65 & 191 & 273\\\hline
\end{tabular}
\end{multicols}
\end{table}

\section{Dependence or independence of parameters}\label{sec:dependence}
Even though the distribution for normalised coefficients up to the bounds we reached computationally can be pretty well described by the GGG function, as seen in the previous section, the result of Radziwi\l\l{} and Soundararajan shows that, for a fixed Hecke eigenform, the parameters must depend on the bound. We investigate this dependence in this section. More precisely, we look at how the parameters behave with respect to~$x$, when we compute with all coefficients up to~$x$.

In order to study this, we consider the case where we have the biggest number of normalised Fourier coefficients, namely $2 \cdot 10^8$, in weight $13/2$. We broke the list of coefficients into 20 subsequent lists of equal size and did the fitting for each of these sublists separately, leading to the results in Table~\ref{T1}.
One clearly observes some dependence. For instance, for the Laplace distribution~\eqref{eq:Laplace}, the $b$ value seems to be slowly, but strictly increasing, whereas the $c$-value slowly, but strictly decreases (with one exception). The values for the Cauchy distribution~\eqref{eq:Cauchy} and the $GG$-distribution~\eqref{eq:GG} also suggest a dependence. For the $GGG$-distribution~\eqref{eq:GGG}, there is a clear dependence of the values for the first couple of sets. However, all four values seem to stabilise for the last sets. The range of data we investigated hence does not allow us to illustrate that the limit distribution is known to be a Dirac delta function.

\begin{table}
\caption{}\label{T1}
\vspace*{0.2cm}
\begin{small}
\begin{tabular}{|l||r|r|r|r||r|r|r||r|r||r|r|r|} 
\hline
&\multicolumn{4}{|c||}{$GGG$} & \multicolumn{3}{|c||}{$GG$} & \multicolumn{2}{|c||}{Laplace} & \multicolumn{3}{|c|}{Cauchy} \\& $a$ & $b$ & $c$ & $d$ & $a$ & $b$ & $c$ & $b$ & $c$ & $a$ & $b$ & $c$ \\\hline\hline$1$ & 0.622 & 1177 & 0.967 & 0.045 & 0.677 & 1038 & 1.08 & 1172 & 0.908 & 142.8 & 0.140 & 0.488 \\\hline$2$ & 0.601 & 1206 & 0.917 & 0.031 & 0.650 & 1082 & 1.009 & 1205 & 0.876 & 142.0 & 0.135 & 0.498 \\\hline$3$ & 0.591 & 1221 & 0.896 & 0.027 & 0.638 & 1101 & 0.981 & 1219 & 0.864 & 143.4 & 0.135 & 0.506 \\\hline$4$ & 0.587 & 1230 & 0.886 & 0.025 & 0.633 & 1110 & 0.969 & 1226 & 0.858 & 144.1 & 0.135 & 0.510 \\\hline$5$ & 0.569 & 1299 & 0.846 & 0.037 & 0.631 & 1117 & 0.960 & 1232 & 0.852 & 144.3 & 0.134 & 0.513 \\\hline$6$ & 0.580 & 1252 & 0.867 & 0.025 & 0.627 & 1123 & 0.952 & 1236 & 0.848 & 144.8 & 0.134 & 0.515 \\\hline$7$ & 0.561 & 1326 & 0.828 & 0.039 & 0.628 & 1126 & 0.948 & 1239 & 0.846 & 144.8 & 0.134 & 0.516 \\\hline$8$ & 0.566 & 1292 & 0.841 & 0.030 & 0.622 & 1133 & 0.940 & 1243 & 0.843 & 145.6 & 0.134 & 0.520 \\\hline$9$ & 0.572 & 1273 & 0.850 & 0.026 & 0.621 & 1136 & 0.936 & 1246 & 0.840 & 145.6 & 0.134 & 0.520 \\\hline$10$ & 0.559 & 1317 & 0.825 & 0.032 & 0.619 & 1141 & 0.930 & 1249 & 0.837 & 146.0 & 0.134 & 0.522\\\hline$11$ & 0.567 & 1289 & 0.839 & 0.027 & 0.619 & 1141 & 0.930 & 1250 & 0.837 & 146.0 & 0.134 & 0.522 \\\hline$12$ & 0.559 & 1317 & 0.824 & 0.031 & 0.618 & 1144 & 0.926 & 1252 & 0.835 & 146.2 & 0.133 & 0.523 \\\hline$13$ & 0.557 & 1326 & 0.819 & 0.032 & 0.617 & 1146 & 0.924 & 1254 & 0.833 & 146.4 & 0.133 & 0.524 \\\hline$14$ & 0.566 & 1275 & 0.840 & 0.020 & 0.610 & 1154 & 0.915 & 1257 & 0.830 & 147.3 & 0.134 & 0.524 \\\hline$15$ & 0.551 & 1349 & 0.807 & 0.034 & 0.616 & 1150 & 0.919 & 1257 & 0.830 & 146.6 & 0.133 & 0.526 \\\hline$16$ & 0.555 & 1332 & 0.814 & 0.031 & 0.616 & 1151 & 0.918 & 1258 & 0.828 & 146.5 & 0.133 & 0.526 \\\hline$17$ & 0.555 & 1327 & 0.814 & 0.029 & 0.613 & 1156 & 0.913 & 1261 & 0.828 & 146.9 & 0.133 & 0.528 \\\hline$18$ & 0.560 & 1309 & 0.822 & 0.025 & 0.612 & 1158 & 0.911 & 1262 & 0.826 & 147.0 & 0.133 & 0.529 \\\hline$19$ & 0.553 & 1333 & 0.811 & 0.029 & 0.612 & 1157 & 0.911 & 1262 & 0.827 & 147.1 & 0.133 & 0.529 \\\hline$20$ & 0.555 & 1323 & 0.814 & 0.026 & 0.610 & 1161 & 0.907 & 1264 & 0.825 & 147.3 & 0.133 & 0.530 \\
\hline
\end{tabular}
\end{small}

\vspace*{0.2cm}
Table of the corresponding RMS-values.
\begin{multicols}{2}\noindent 
\begin{tabular}{|l||r|r|r|r|}\hline& $GGG$ & $GG$ & Laplace & Cauchy \\
\hline
\hline$1$ & 18.4 & 18.9 & 39.0 & 33.6\\\hline$2$ & 18.1 & 18.6 & 35.4 & 31.6\\\hline$3$ & 18.7 & 19.2 & 34.1 & 31.2\\\hline$4$ & 18.2 & 18.8 & 33.3 & 30.7\\\hline$5$ & 18.0 & 18.8 & 33.2 & 30.0\\\hline$6$ & 18.5 & 19.0 & 32.7 & 30.4\\\hline$7$ & 18.4 & 19.3 & 33.0 & 29.9\\\hline$8$ & 18.2 & 19.0 & 31.9 & 29.7\\\hline$9$ & 18.3 & 18.8 & 31.7 & 29.9\\\hline$10$ & 18.3 & 19.1 & 31.6 & 29.5\\\hline
\end{tabular}

\noindent 
\begin{tabular}{|l||r|r|r|r|}
\hline& 
$GGG$ & $GG$ & Laplace & Cauchy \\
\hline
\hline$11$ & 18.4 & 19.0 & 31.6 & 29.8\\\hline$12$ & 18.3 & 19.2 & 31.5 & 29.5\\\hline$13$ & 18.3 & 19.2 & 31.5 & 29.5\\\hline$14$ & 18.2 & 18.7 & 30.1 & 29.4\\\hline$15$ & 17.7 & 18.7 & 31.0 & 28.9\\\hline$16$ & 18.3 & 19.2 & 31.3 & 29.4\\\hline$17$ & 18.0 & 18.8 & 30.6 & 29.1\\\hline$18$ & 18.3 & 19.0 & 30.7 & 29.4\\\hline$19$ & 18.4 & 19.2 & 30.8 & 29.2\\\hline$20$ & 18.4 & 19.2 & 30.5 & 29.2\\\hline\end{tabular}
\end{multicols}
\end{table}

\begin{figure}[h]
\centering
\includegraphics[width=9cm]{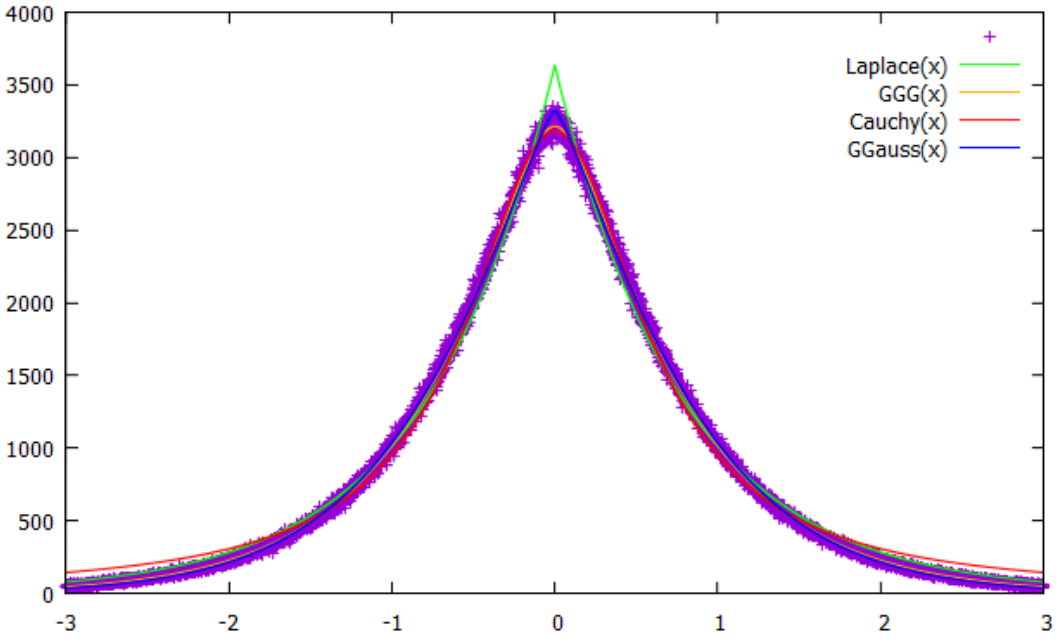}
\caption{Histogram and distributions of Hecke eigenform of weight $13/2$ for prime indexed coefficients only}
\label{fig:p13}
\end{figure}

Recall that we restricted our efforts to squarefree indices $n$ such that the coefficient is not automatically known to be zero. We investigated further if the distribution seems to change significantly when considering prime indexed coefficients only. In the integral weight case, there are huge differences and the semi-circular distribution of Sato-Tate for non-CM eigenforms is only valid for prime indices. Moreover, the coefficients of integral weight Hecke eigenforms are multiplicative functions, hence the distribution of the coefficients at prime indices determines the rest. We do not know of any reason to believe that similar things happen in the half-integral weight case. Indeed, the shapes of the distribution graphs do not change significantly if we restrict to prime indices (see Figure~\ref{fig:p13} for an example). Of course, some of the best fit parameters move, as we can see in Table~\ref{T4}.

\begin{table}[H]
\caption{}\label{T4}
\vspace*{0.2cm}
\begin{small}
\begin{tabular}{|l||l|l|l|l||l|l|l||l|l||l|l|l|}
\hline
& \multicolumn{4}{|c||}{$GGG$} & \multicolumn{3}{|c||}{$GG$} & \multicolumn{2}{|c||}{Laplace} & \multicolumn{3}{|c|}{Cauchy} \\& $a$ & $b$ & $c$ & $d$ & $a$ & $b$ & $c$ & $b$ & $c$ & $a$ & $b$ & $c$ \\
\hline
\hline
Sqfree & 0.570 & 25666 & 0.850 & 0.030 & 0.623 & 22621 & 0.942 & 24837 & 0.843 & 140 & 0.006 & 0.113 \\
\hline
Prime & 0.550 & 3892 & 0.759 & 0.030 & 0.614 & 3329 & 0.856 & 3635 & 0.784 & 196 & 0.062 & 0.380 \\
\hline
\end{tabular}
\end{small}
\end{table}

The $b$ values differ just because there are far more squarefree numbers than prime numbers. For the same reason also the RMS values are different. However, the most important parameters, i.e.\ the $a$-parameters in $GGG$ and $GG$ are very similar. It is not clear if the slight change of parameters can be explained by the fact that the set of primes among squarefree numbers is biassed towards small values, or if the change is not significant. This gives us confidence in our belief that prime indexed normalised coefficients are not distributed differently from those with squarefree indices, when considering coefficients with indices up to the bounds we could reach.

The observed very good approximations by the $GGG$-distribution of the normalised Fourier coefficients up to the bounds we could computationally reach lead us to ask the question whether this holds for all bounds.
\newpage

\begin{question}\label{qu:GGG} Let $f$ be a half-integral weight cuspidal Hecke eigenform in the Kohnen plus-space in half-integral weight $k=\ell + \frac{1}{2}$ and let $b(n)$ be its normalised coefficients. Let $x \in \mathbb{R}_{>0}$.
Can the distribution of the $b(n)$ for $n \le x$ squarefree and $n \equiv (-1)^\ell \mod 4$ be approximated by the $GGG$-distribution with parameters depending on $f$ and~$x$?

More precisely, are there constants $a,c,d$, depending on $f$ and~$x$, such that for all intervals $I=[\alpha,\beta] \subseteq \mathbb{R}$ the quotient\[ \frac{\#\{n \in \mathbb{R} \textnormal{ squarefree }|\; n \le x, n \equiv (-1)^\ell \mod 4, b(n) \in I \}}{\#\{n \in \mathbb{N} \textnormal{ squarefree }|\; n \le x, n \equiv (-1)^\ell \mod 4\}}\]is `close to'\[ \frac{1}{b} \int_\alpha^\beta \exp(-\frac{(d+t^2)^a}{c}) dt,\]where $b = \int_{-\infty}^\infty \exp(-\frac{(d+t^2)^a}{c}) dt$?
\end{question}

\section{A strengthening of the Bruinier-Kohnen conjecture}\label{sec:BK}
Bruinier and Kohnen conjectured that the signs of coefficients of half-integral weight Hecke eigenforms are equidistributed. More precisely, let $f=\sum_{n=1}^{\infty} a(n) q^n \in S_k(N, \chi)$ be a cusp form of weight $k=\ell+1/2$ with real Fourier coefficients and assume that $f$ is orthogonal to the unary theta series when $\ell=1$. Then the Bruinier-Kohnen Conjecture (\cite{10} and \cite{28}) asserts that the sets $\{ n \in \mathbb{N} : a(n) > 0 \}$ and $\{ n \in \mathbb{N} : a(n) < 0 \}$ have the same natural density, equal to half of the natural density of $\{ n \in \mathbb{N} : a(n) \neq 0 \}$. 

Combining the Shimura lift with the (proved) celebrated Sato-Tate Conjecture for integral weight Hecke eigenforms, it is not very difficult to prove equidistribution of signs for the coefficients indexed by squares, see \cite{6}, \cite{22}, \cite{23}. The sign equidistribution problem has still received much attention and is widely studied (for instance \cite{7}), and the technique from \cite{22} has been extended to more general automorphic forms like Hilbert modular forms in \cite{26}. Note that this is only a partial result and the full proof of the conjecture is still an open problem and, for the moment, it is likely that there is no theoretical tool to attack this problem.

We see the calculations in this article as the most systematic and largest computational support for the Bruinier-Kohnen Conjecture to date.In fact, if the distribution of coefficients up to any bound $x$ follows any distribution function (depending on $x$) that is symmetric with respect to~$0$, e.g.any of the four types discussed above, then the Bruinier-Kohnen Conjecture is true.

In fact, the symmetry around $0$ suggests that the signs are uniformly distributed and that the distribution of the signs and the distribution of the absolute value of the normalised coefficients are independent. In order to make this precise, we recall that according to \eqref{eq:lower}, there are infinitely many `big' normalised coefficients $|b(n)|$. This suggests that any non-empty interval $I \subseteq \mathbb{R}_{> 0}$ will contain infinitely many normalised coefficients $|b(n)|$ for squarefree~$n$. We insist on squarefree because we do not want to deal with the contributions that are understood by the Shimura lift.

We feel that the symmetry around $0$ of the distribution of the normalised coefficients up to the bounds that we computationally reached warrant the formulation of the following conjecture, strengthening the Bruinier-Kohnen Conjecture.

\begin{conjecture}[Independence of sign and absolute value] 
\label{conj:SBK} Let $f$ be as above. Let $I \subseteq \mathbb{R}_{> 0}$ be any interval.Then the following limit exists and we have\[ \lim_{x \to \infty} \frac{\#\{ n \le x \textnormal{ squarefree } \;|\; |b(n)| \in I, b(n) > 0\} }{\#\{ n \le x \textnormal{ squarefree } \;|\; |b(n)| \in I\} } = \frac{1}{2}.\]\end{conjecture}

\section*{Acknowledgment}This work is supported by The Scientific and Technological Research Council of Turkey (TUBITAK) with the project number 118F148.I.I.\ acknowledges partial and complement support by Bilecik Seyh Edebali University research project number 2018-01.BSEU.04-01 and would like to thank the University of Luxembourg for the great hospitality in several visits.The authors thank the Izmir Institute of Technology and especially Prof.\ Dr.\ Engin B\"{u}y\"{u}ka\c{s}{\i}k, where a huge part of this research was carried out, for its great hospitality.
The authors are extremely grateful to Kannan Soundararajan for very useful feedback on a first version of this article, clarifying various points to them. They also thank Henri Cohen for providing them with the initial code for obtaining Hecke eigenforms via Rankin-Cohen brackets in Pari/GP and Winfried Kohnen for interesting discussions. They also thank Oktay Pashaev for some interesting suggestions. Thanks are also due to the anonymous referees for helpful suggestions concerning the presentation of the paper. I.I.\ wishes to thank his parents who welcomed him during his Izmir visit in summer 2019. G.W.\ thanks Giovanni Peccati for suggesting the generalised Gaussian distribution.

\newpage\section*{Appendix: Graphs of histograms of all computed examples}
\begin{figure}[H]\centering\vspace*{-2.6cm}\hspace*{-2cm}\includegraphics[width=1.1\textwidth]{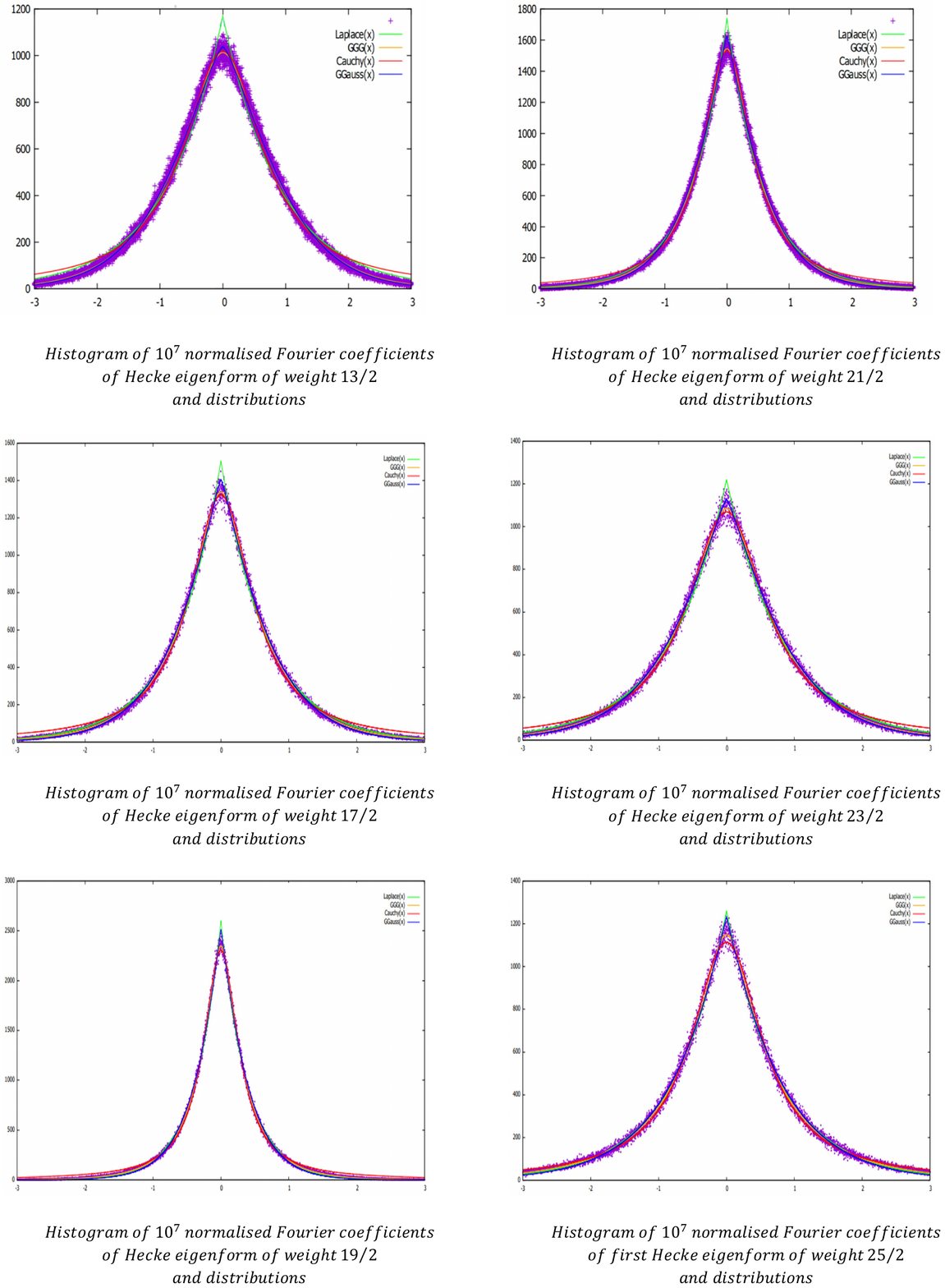}
\end{figure}
\begin{figure}[H]\centering\vspace*{-2.6cm}\hspace*{-2cm}\includegraphics[width=1.1\textwidth]{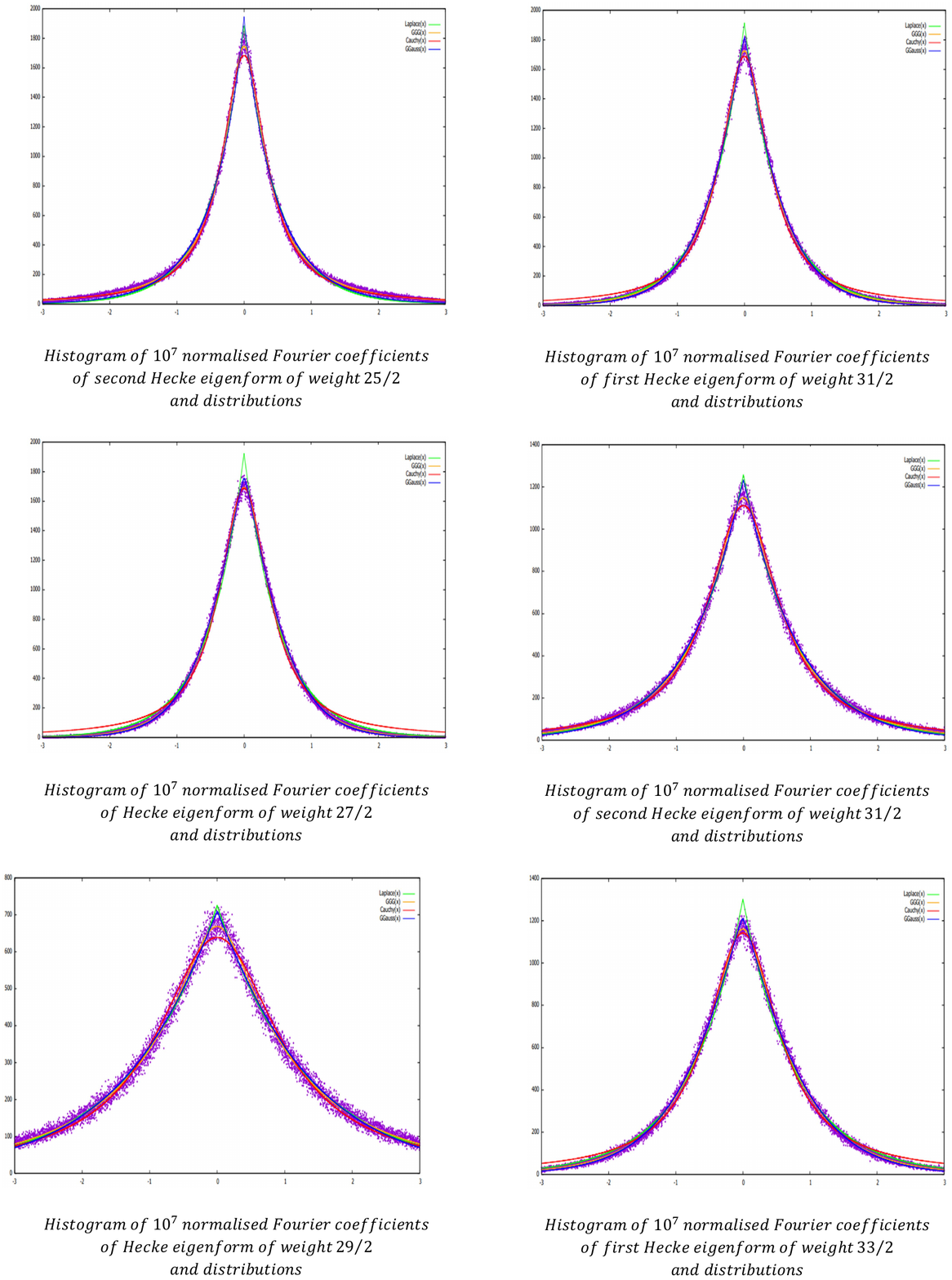}
\end{figure}
\begin{figure}[H]\centering\vspace*{-2.6cm}\hspace*{-2cm}\includegraphics[width=1.1\textwidth]{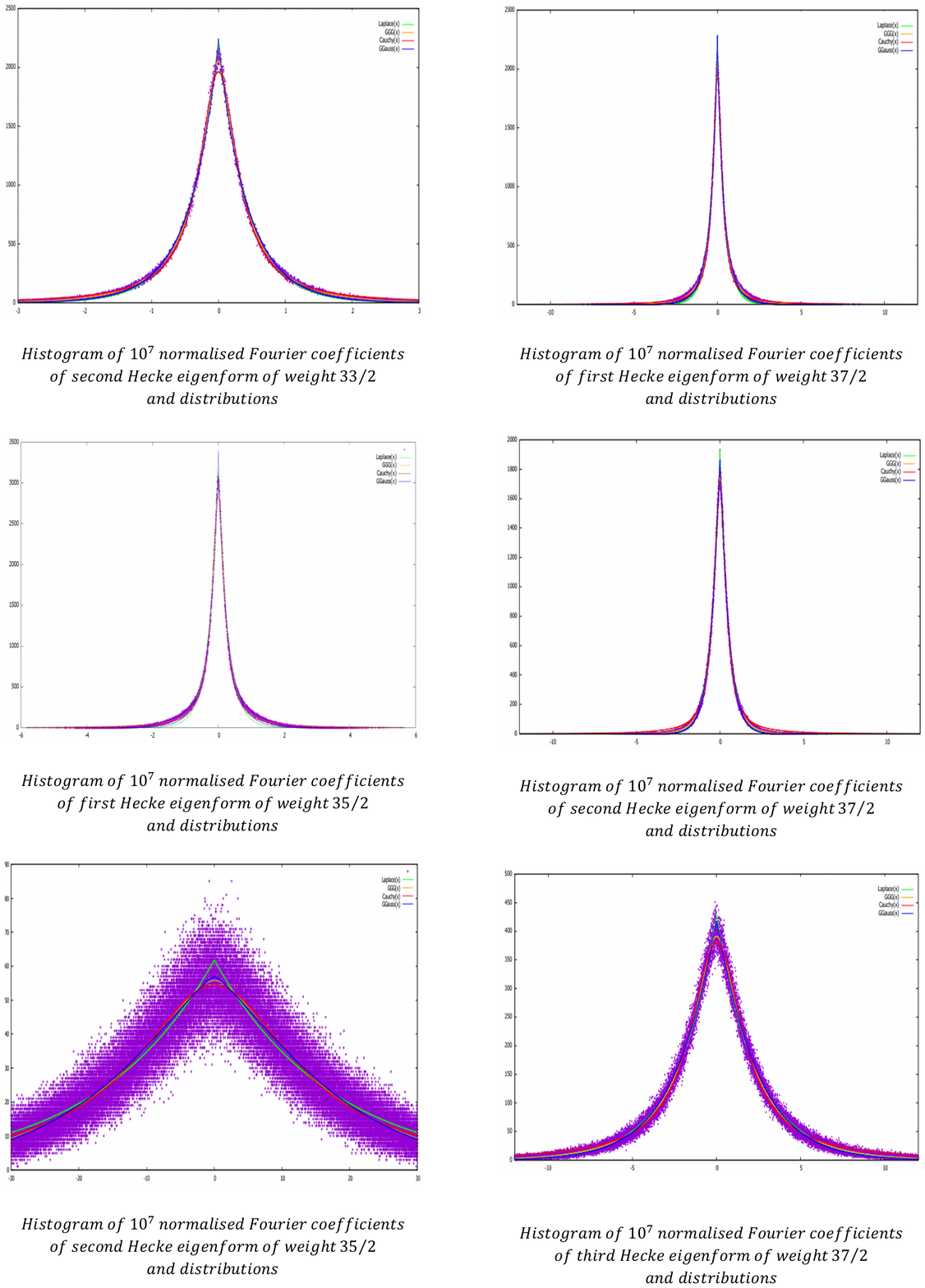}
\end{figure}
\begin{figure}[H]\centering\vspace*{-2.6cm}\hspace*{-2cm}\includegraphics[width=1.1\textwidth]{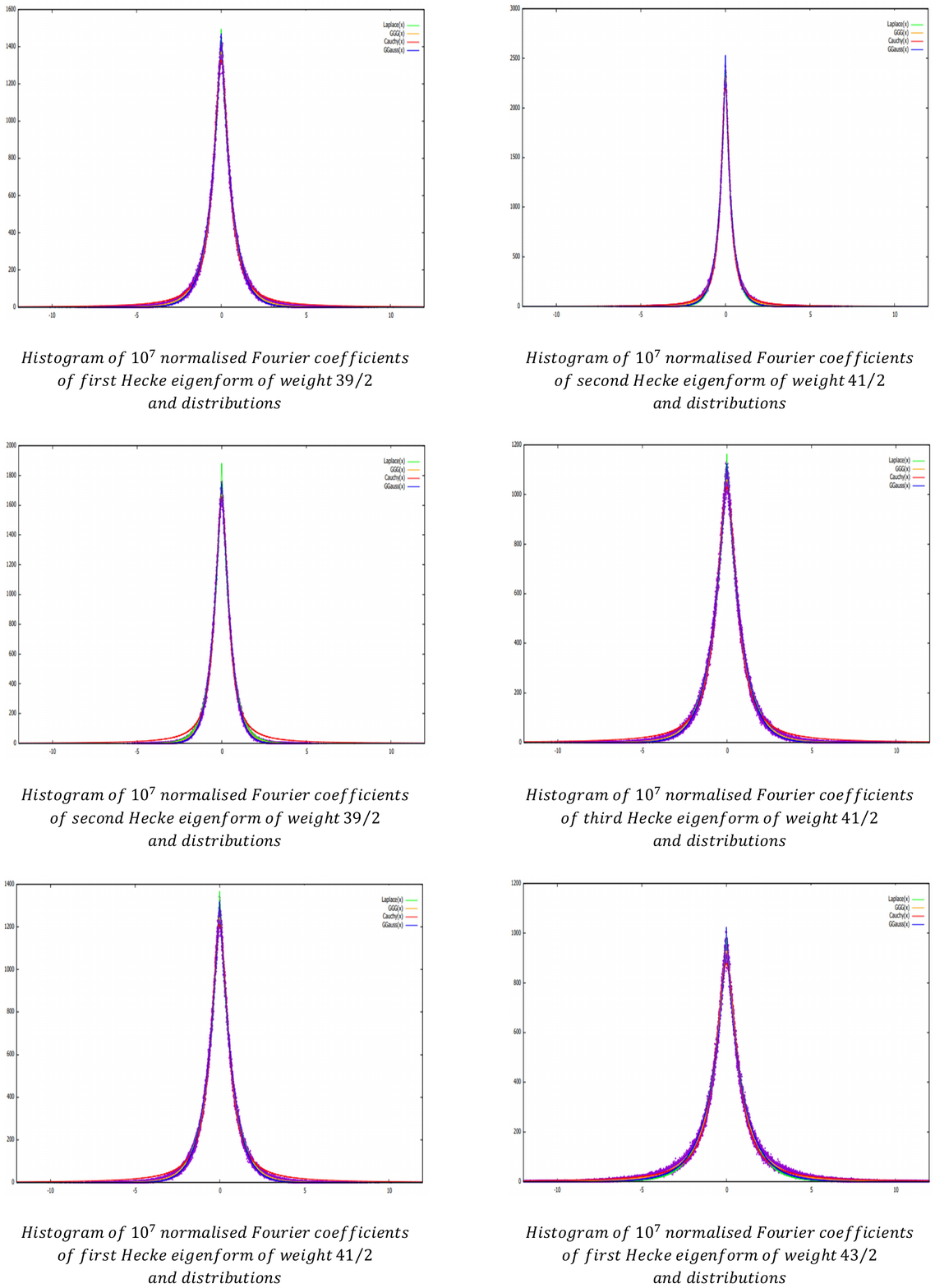}
\end{figure}
\begin{figure}[H]\centering\vspace*{-2.6cm}\hspace*{-2cm}\includegraphics[width=1.1\textwidth]{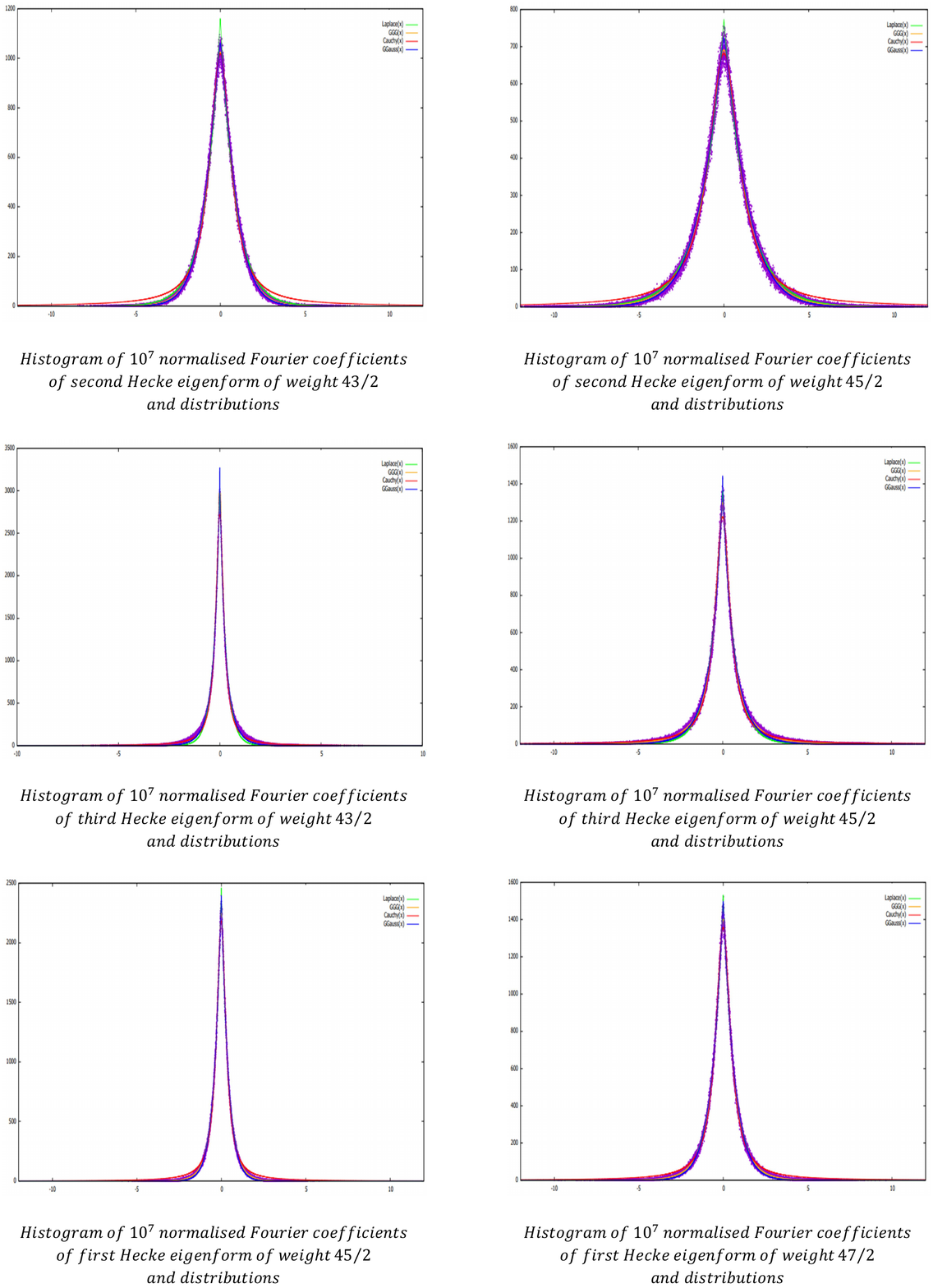}
\end{figure}
\begin{figure}[H]\centering\vspace*{-2.6cm}\hspace*{-2cm}\includegraphics[width=1.1\textwidth]{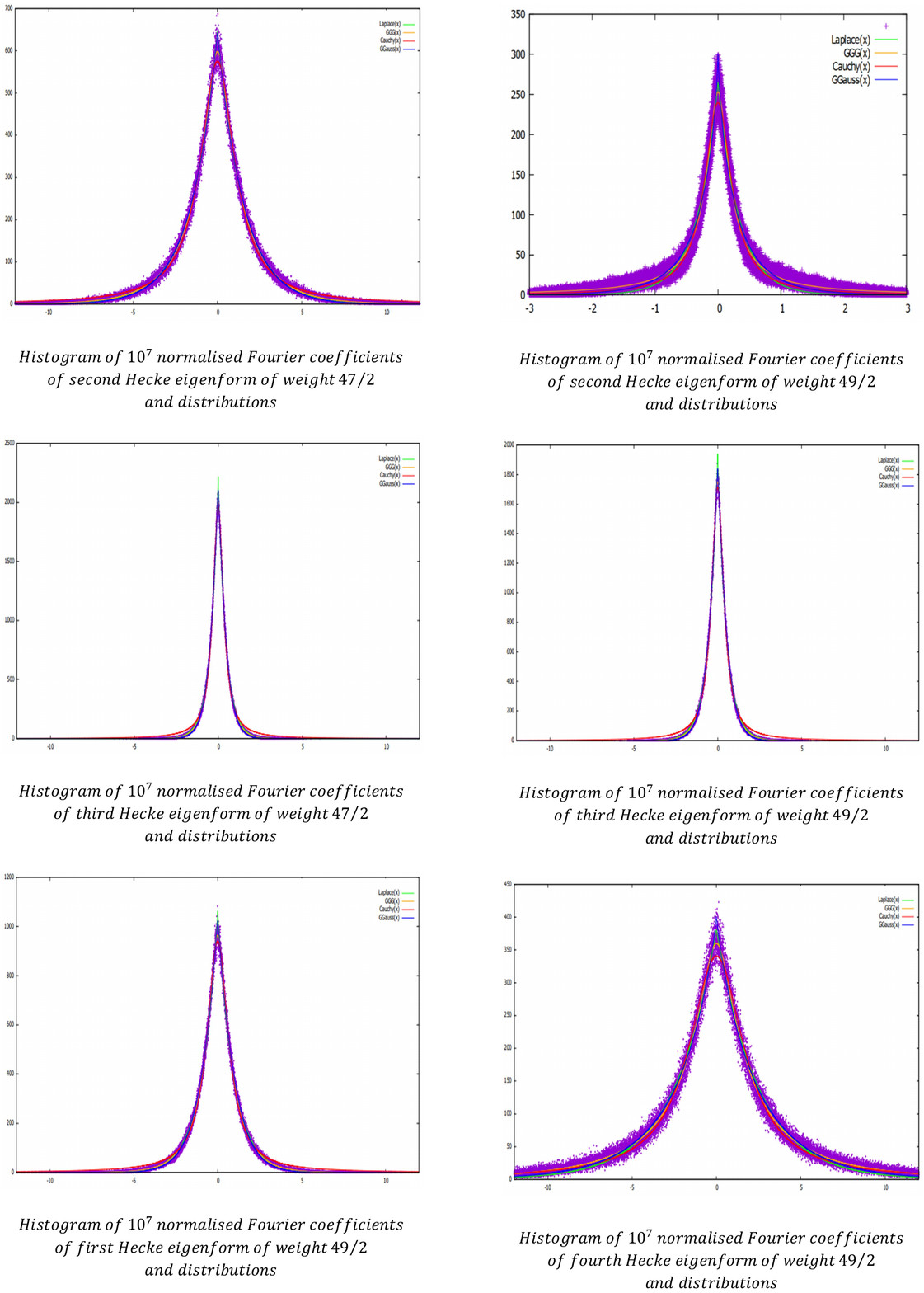}
\end{figure}
\begin{figure}[H]\centering\vspace*{-2.6cm}\hspace*{-2cm}\includegraphics[width=1.1\textwidth]{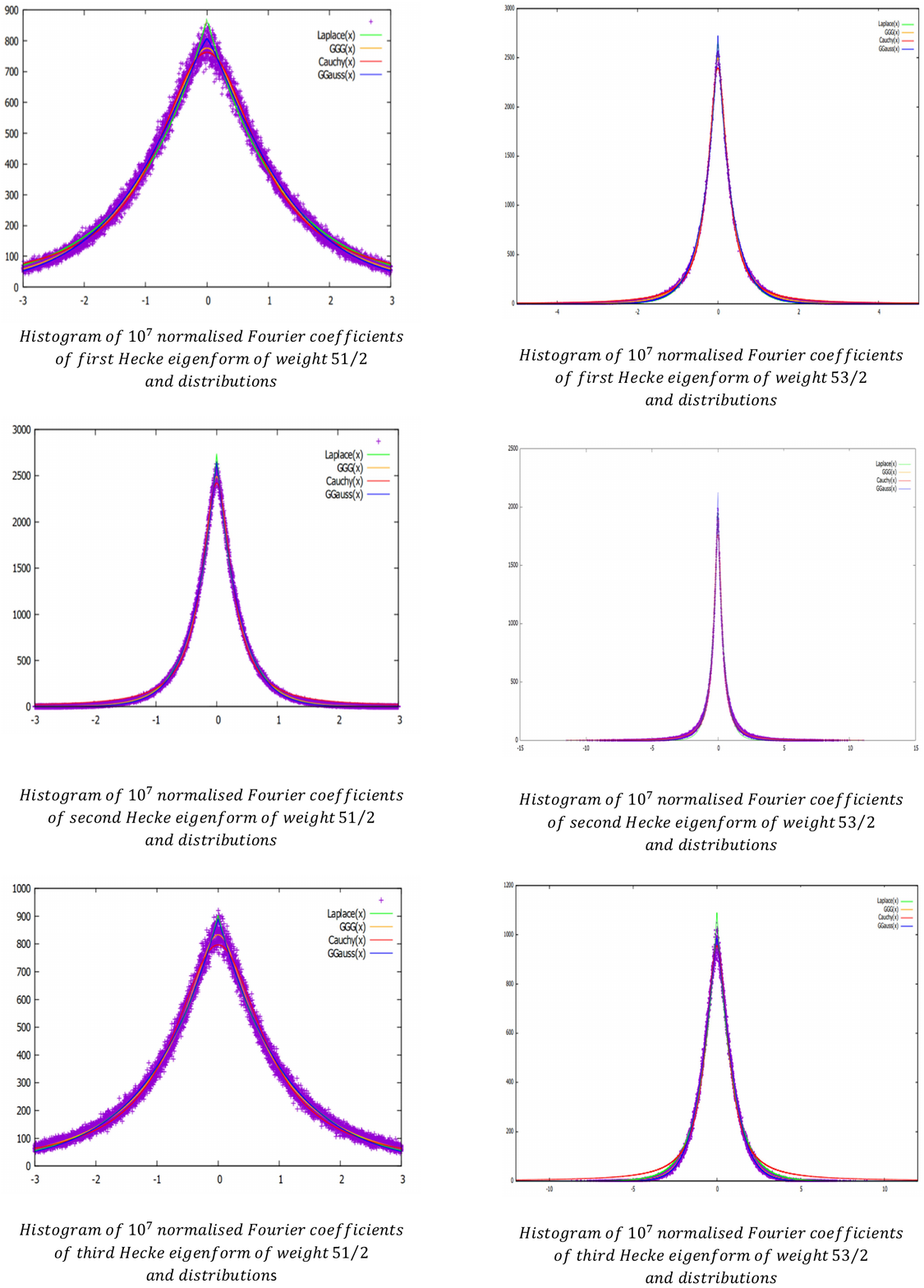}
\end{figure}
\begin{figure}[H]\centering\vspace*{-2.6cm}\hspace*{-2cm}\includegraphics[width=1.1\textwidth]{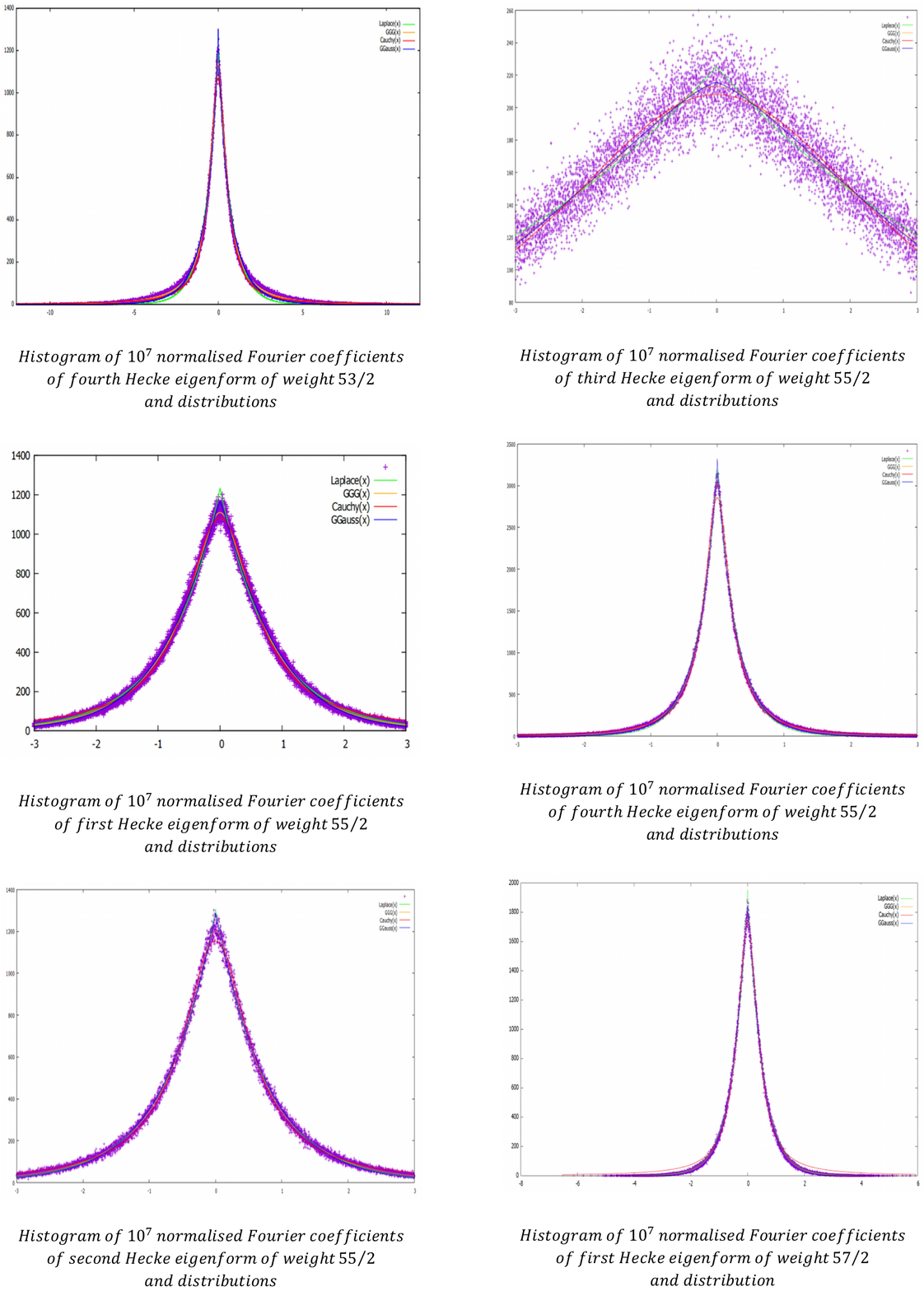}
\end{figure}
\begin{figure}[H]\centering\vspace*{-2.6cm}\hspace*{-2cm}\includegraphics[width=1.1\textwidth]{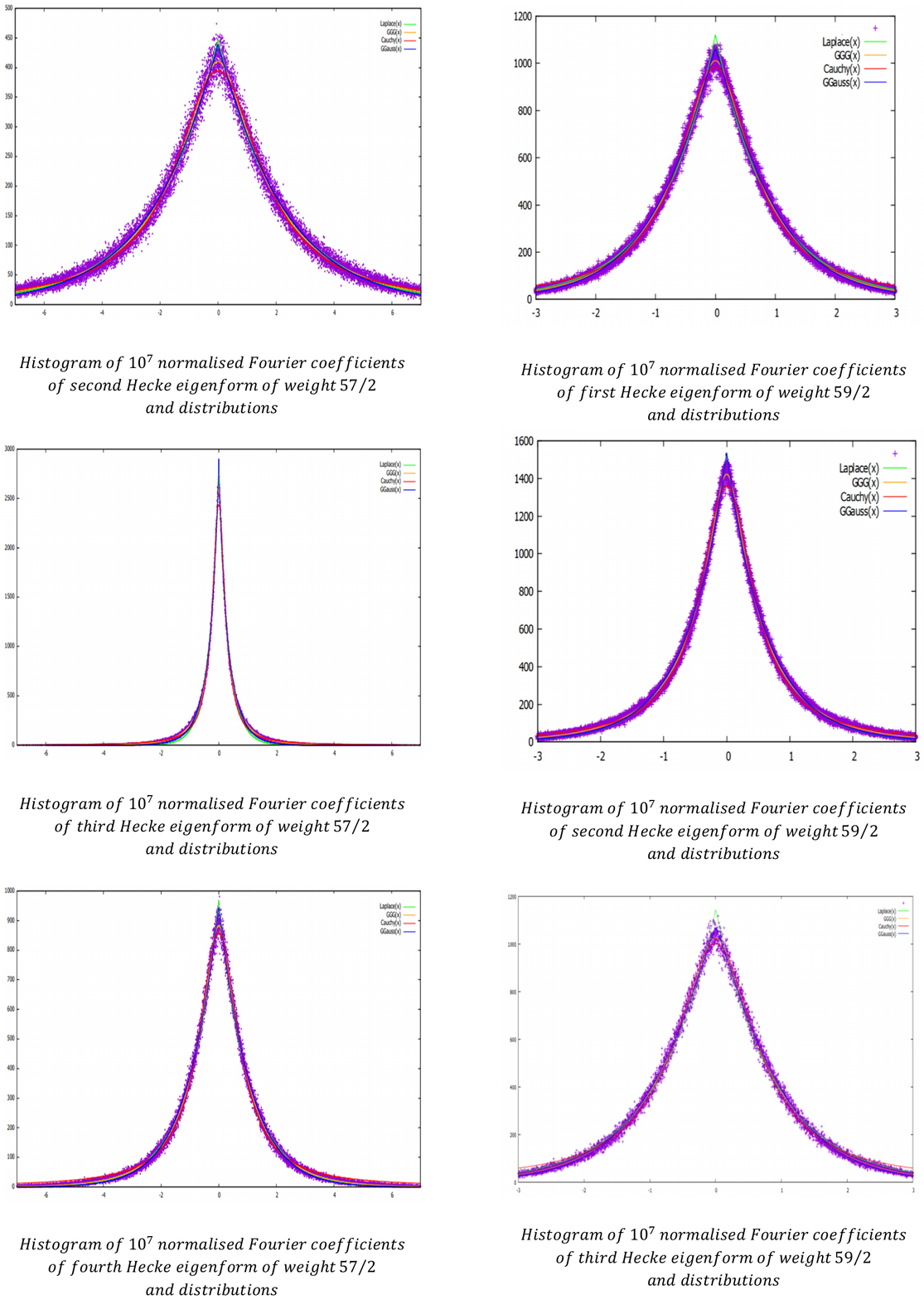}
\end{figure}
\begin{figure}[H]\centering\vspace*{-2.6cm}\hspace*{-2cm}\includegraphics[width=1.1\textwidth]{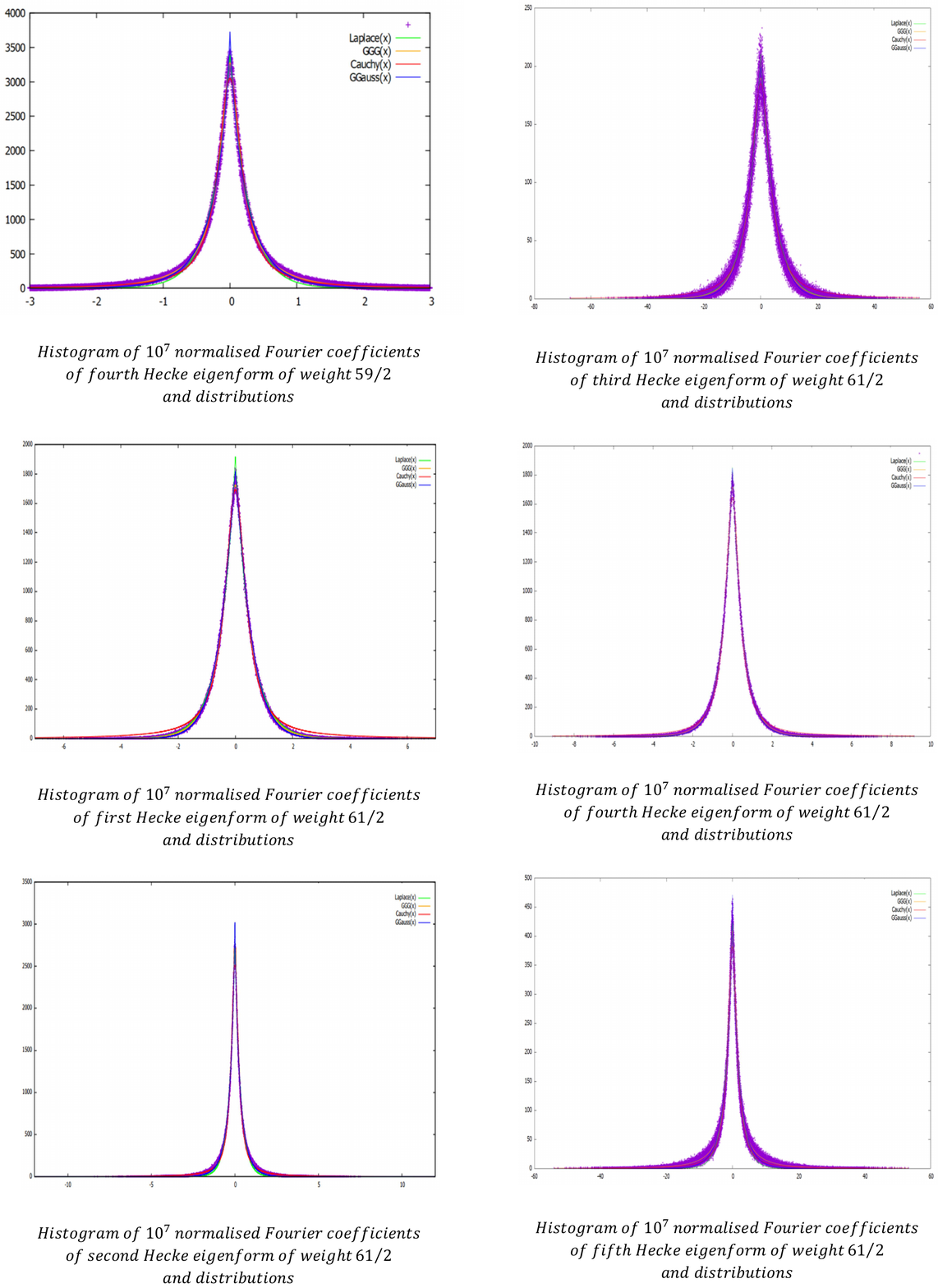}
\end{figure}
\begin{figure}[H]\centering\vspace*{-2.6cm}\hspace*{-2cm}\includegraphics[width=1.1\textwidth]{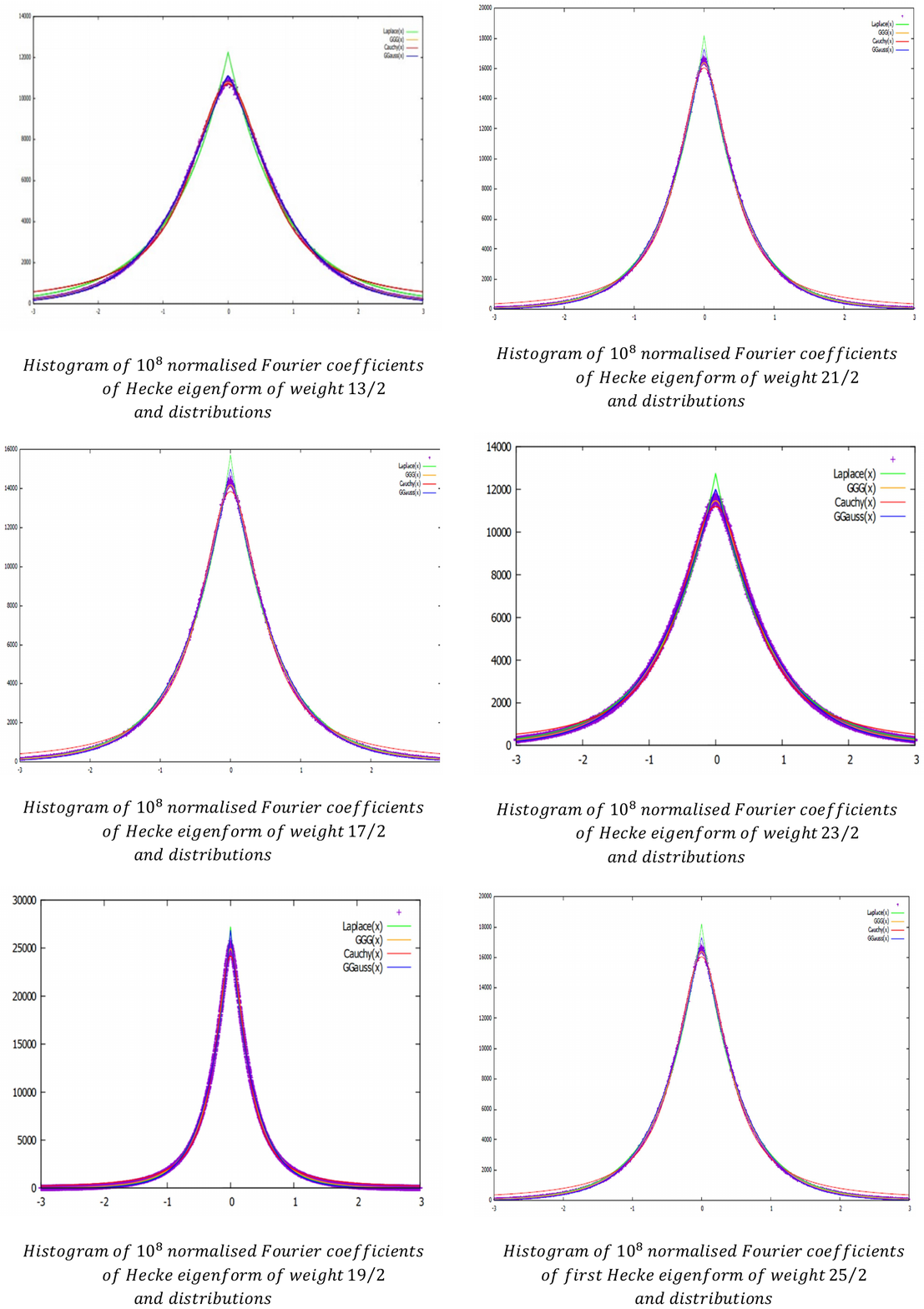}
\end{figure}
\begin{figure}[H]\centering\vspace*{-2.6cm}\hspace*{-2cm}\includegraphics[width=1.1\textwidth]{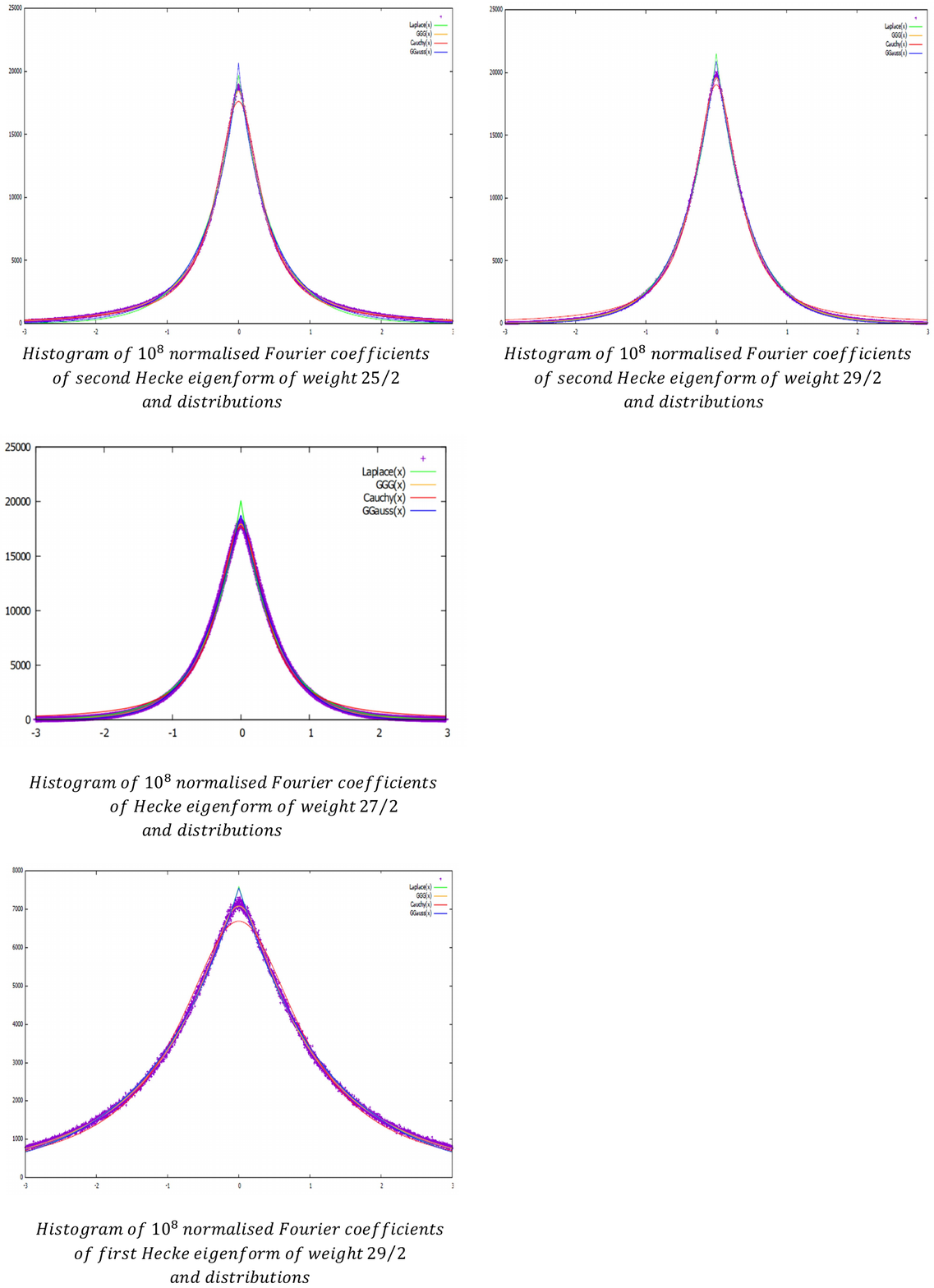}
\end{figure}
\newpage\section*{Appendix: Graphs of histograms in weight $13/2$ with $2 \cdot 10^8$ coefficients in $20$ subsets}
\begin{figure}[H]\centering\vspace*{-2.6cm}\hspace*{-2cm}\includegraphics[width=1.1\textwidth]{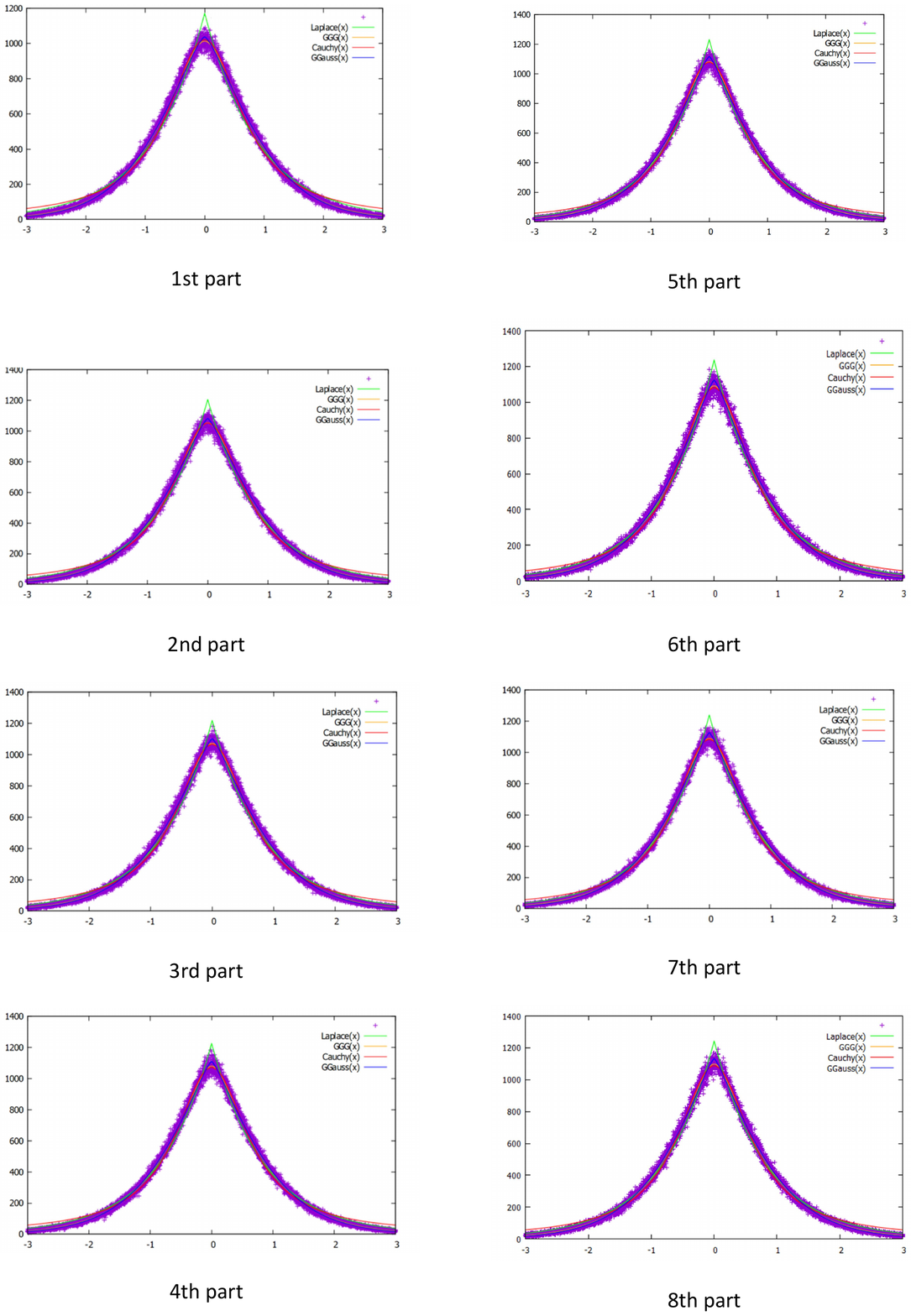}
\end{figure}
\begin{figure}[H]\centering\vspace*{-2.6cm}\hspace*{-2cm}\includegraphics[width=1.1\textwidth]{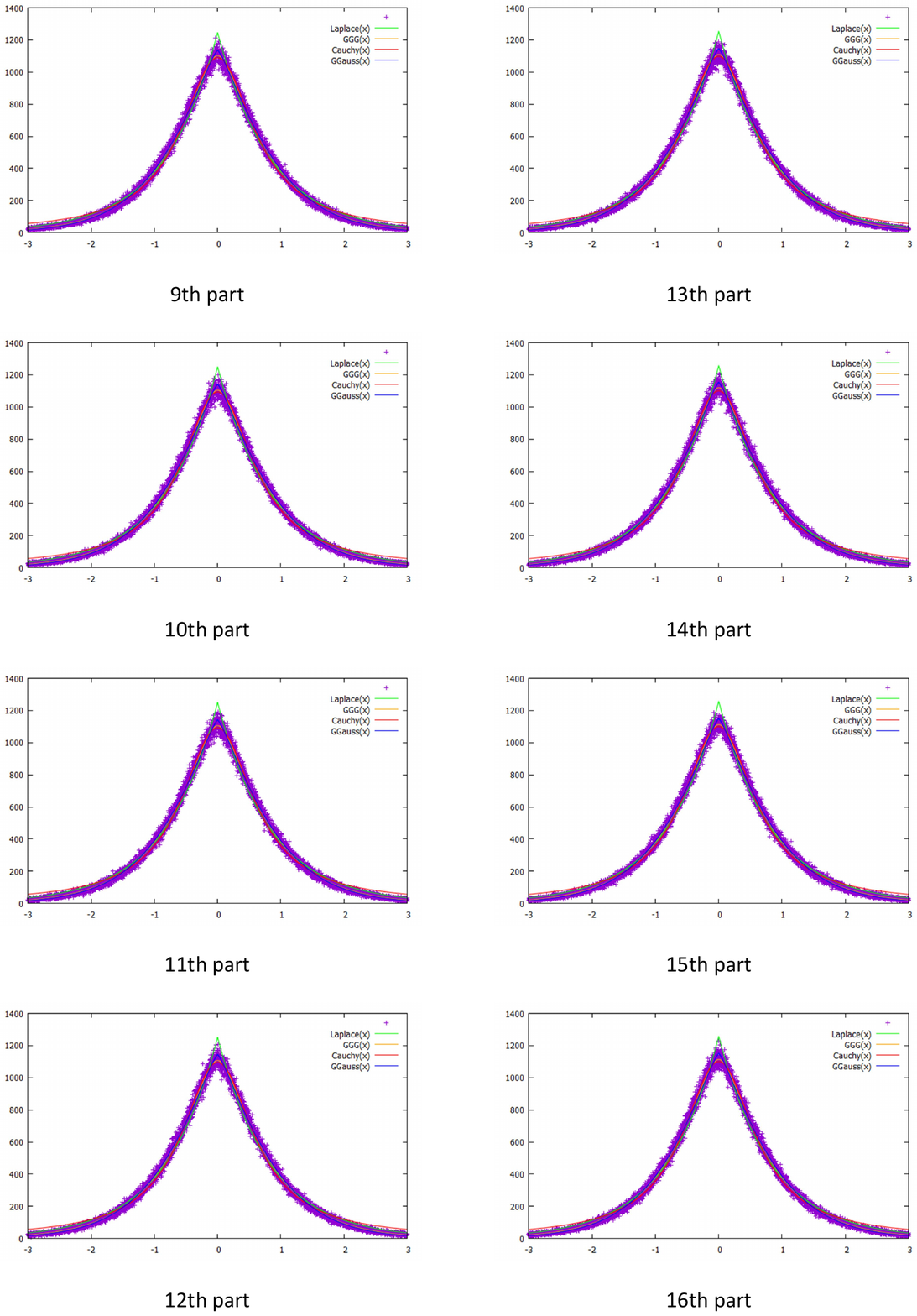}
\end{figure}
\begin{figure}[H]\centering\vspace*{-2.6cm}\hspace*{-2cm}\includegraphics[width=1.1\textwidth]{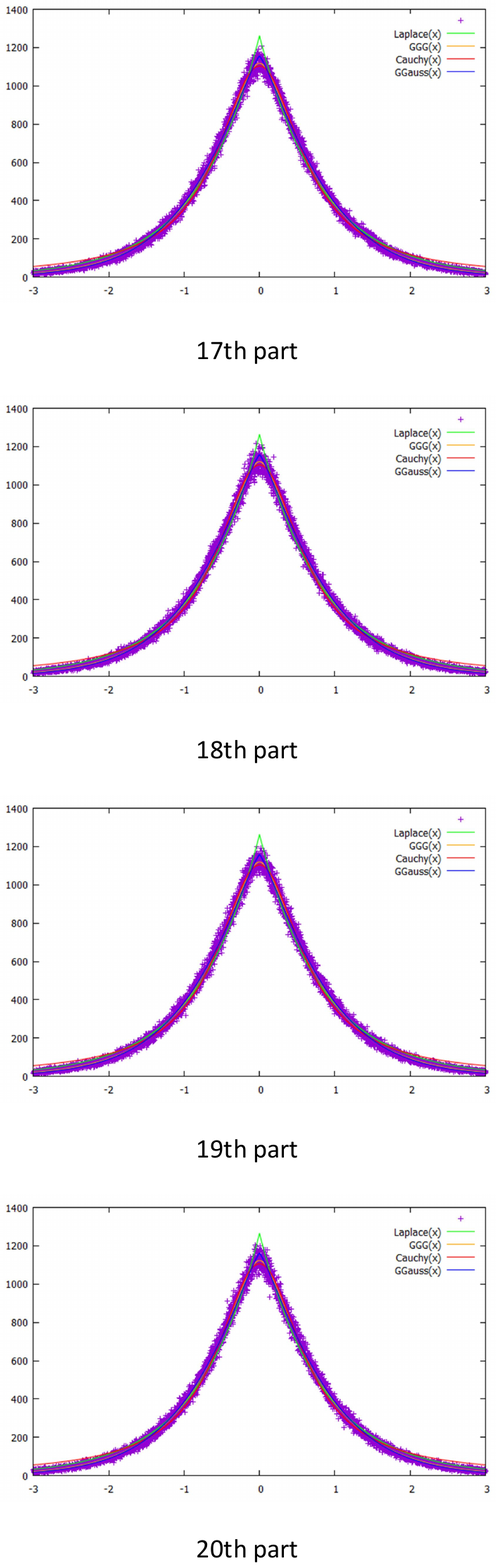}
\end{figure}

\newpage\section*{Appendix: Tables of best fit parameters}
\noindent Best fit parameters (rounded) for the $GGG$-distribution for all examples with $10^7$ coefficients:
\begin{multicols}{2}
\begin{small}
\noindent
\begin{tabular}{|l||r|r|r|r|}\hline& $a$ & $b$ & $c$ & $d$\\\hline\hline$13/2$ & 0.622 & 1177.4 & 0.967 & 0.045\\\hline$17/2$ & 0.470 & 1986.6 & 0.575 & 0.043\\\hline$19/2$ & 0.386 & 4595.4 & 0.318 & 0.018\\\hline$21/2$ & 0.477 & 2268.8 & 0.506 & 0.032\\\hline$23/2$ & 0.527 & 1358.6 & 0.800 & 0.039\\\hline$25/2(1)$ & 0.384 & 2065.3 & 0.570 & 0.057\\\hline$25/2(2)$ & 0.219 & 12428.5 & 0.262 & 0.048\\\hline$27/2$ & 0.542 & 2129.2 & 0.498 & 0.018\\\hline$29/2$ & 0.354 & 1272.2 & 0.789 & 0.147\\\hline$31/2(1)$ & 0.468 & 2404.4 & 0.469 & 0.019\\\hline$31/2(2)$ & 0.375 & 2162.2 & 0.555 & 0.061\\\hline$33/2(1)$ & 0.508 & 1510.6 & 0.721 & 0.038\\\hline$33/2(2)$ & 0.338 & 4206.6 & 0.3526 & 0.017\\\hline$35/2(1)$ & 0.185 & 30546.8 & 0.195 & 0.014\\\hline$35/2(2)$ & 0.595 & 61.6 & 30.3975 & 6.408 \\\hline$37/2(1)$ & 0.248 & 7668.0 & 0.292 & 0.022\\\hline$37/2(2)$ & 0.384 & 3432.2 & 0.397 & 0.035\\\hline$37/2(3)$ & 0.414 & 620.5 & 1.507 & 0.415\\\hline$39/2(1)$ & 0.397 & 2286.9 & 0.519 & 0.035\\\hline$39/2(2)$ & 0.508 & 2217.3 & 0.493 & 0.021\\\hline$41/2(1)$ & 0.439 & 1830.7 & 0.609 & 0.037\\\hline$41/2(2)$ & 0.334 & 4608.3 & 0.329 & 0.012\\\hline$41/2(3)$ & 0.441 & 1534.6 & 0.708 & 0.048\\\hline$43/2(1)$ & 0.307 & 2131.8 & 0.552 & 0.080\\\hline$43/2(2)$ & 0.548 & 1252.9 & 0.879 & 0.043\\\hline$43/2(3)$ & 0.232 & 12488.3 & 0.238 & 0.011\\\hline$45/2(1)$ & 0.419 & 3505.8 & 0.357 & 0.012\\\hline$45/2(2)$ & 0.492 & 928.6 & 1.172 & 0.116\\\hline$45/2(3)$ & 0.299 & 3046.2 & 0.443 & 0.038\\\hline$47/2(1)$ & 0.408 & 2253.6 & 0.521 & 0.033\\\hline
\end{tabular}

\noindent
\begin{tabular}{|l||r|r|r|r|}\hline& $a$ & $b$ & $c$ & $d$\\\hline\hline$47/2(2)$ & 0.403 & 929.1 & 1.049 & 0.147\\\hline$47/2(3)$ & 0.475 & 2714.2 & 0.412 & 0.013\\\hline$49/2(1)$ & 0.439 & 1428.2 & 0.760 & 0.064\\\hline$49/2(2)$ & 0.094 & 798896 & 0.121 & 0.022\\\hline$49/2(3)$ & 0.480 & 2317.8 & 0.474 & 0.015\\\hline$49/2(4)$ & 0.269 & 1124.0 & 0.782 & 0.648\\\hline$51/2(1)$ & 0.509 & 1026.3 & 1.075 & 0.094\\\hline$51/2(2)$ & 0.442 & 3572.5 & 0.335 & 0.008 \\\hline$51/2(3)$ & 0.369 & 1488.4 & 0.713 & 0.091\\\hline$53/2(1)$ & 0.339 & 5188.6 & 0.306 & 0.012\\\hline$53/2(2)$ & 0.274 & 4932.5 & 0.341 & 0.017\\\hline$53/2(3)$ & 0.570 & 1143.5 & 0.979 & 0.0478\\\hline$53/2(4)$ & 0.220 & 5509.5 & 0.350 & 0.067\\\hline$55/2(1)$ & 0.475 & 1518.7 & 0.718 & 0.043\\\hline$55/2(2)$ & 0.392 & 1998.2 & 0.574 & 0.043\\\hline$55/2(3)$ & 0.567 & 235.1 & 5.046 & 0.286\\\hline$55/2(4)$ & 0.338 & 5684.9 & 0.280 & 0.006\\\hline$57/2(1)$ & 0.488 & 2320.4 & 0.473 & 0.016\\\hline$57/2(2)$ & 0.378 & 735.0 & 1.237 & 0.422\\\hline$57/2(3)$ & 0.276 & 7174.6 & 0.281 & 0.011\\\hline$57/2(4)$ & 0.415 & 1401.3 & 0.771 & 0.083\\\hline$59/2(1)$ & 0.479 & 1351.5 & 0.796 & 0.047\\\hline$59/2(2)$ & 0.365 & 2562.3 & 0.481 & 0.031\\\hline$59/2(3)$ & 0.515 & 1281.0 & 0.841 & 0.038\\\hline$59/2(4)$ & 0.249 & 10585.4 & 0.241 & 0.006\\\hline$61/2(1)$ & 0.465 & 2445.7 & 0.465 & 0.020\\\hline$61/2(2)$ & 0.283 & 6542.5 & 0.286 & 0.008\\\hline$61/2(3)$ & 0.450 & 285.4 & 3.280 & 1.910\\\hline$61/2(4)$ & 0.395 & 2718.4 & 0.445 & 0.018\\\hline$61/2(5)$ & 0.167 & 4656.1 & 0.375 & 0.580\\\hline
\end{tabular}
\end{small}
\end{multicols}
\newpage\noindent Best fit parameters (rounded) for the $GG$-distribution for all examples with $10^7$ coefficients:\begin{multicols}{2}\noindent\begin{tabular}{|l||r|r|r|}\hline& $a$ & $b$ & $c$\\\hline\hline$13/2$ & 0.677 & 1038 & 1.08\\\hline$17/2$ & 0.581 & 1406 & 0.71\\\hline$19/2$ & 0.538 & 2513 & 0.37\\\hline$21/2$ & 0.585 & 1622 & 0.60\\\hline$23/2$ & 0.599 & 1128 & 0.93\\\hline$25/2(1)$ & 0.525 & 1232 & 0.80\\\hline$25/2(2)$ & 0.467 & 1946 & 0.51\\\hline$27/2$ & 0.615 & 1756 & 0.54\\\hline$29/2$ & 0.529 & 708 & 1.37\\\hline$31/2(1)$ & 0.560 & 1817 & 0.53\\\hline$31/2(2)$ & 0.523 & 1232 & 0.80\\\hline$33/2(1)$ & 0.590 & 1211 & 0.85\\\hline$33/2(2)$ & 0.481 & 2243 & 0.44\\\hline$35/2(1)$ & 0.428 & 3402 & 0.32\\\hline$35/2(2)$ & 0.646 & 57 & 43.86 \\\hline$37/2(1)$ & 0.437 & 2291 & 0.45\\\hline$37/2(2)$ & 0.540 & 1862 & 0.52\\\hline$37/2(3)$ & 0.529 & 418 & 2.55\\\hline$39/2(1)$ & 0.518 & 1470 & 0.68\\\hline$39/2(2)$ & 0.593 & 1744 & 0.56\\\hline$41/2(1)$ & 0.540 & 1319 & 0.76\\\hline$41/2(2)$ & 0.471 & 2532 & 0.40\\\hline$41/2(3)$ & 0.539 & 1123 & 0.91\\\hline$43/2(1)$ & 0.458 & 1025 & 0.93\\\hline$43/2(2)$ & 0.614 & 1061 & 1.02\\\hline$43/2(3)$ & 0.426 & 3277 & 0.33\\\hline$45/2(1)$ & 0.528 & 2396 & 0.40\\\hline$45/2(2)$ & 0.579 & 724 & 1.55\\\hline$45/2(3)$ & 0.448 & 1444 & 0.67\\\hline$47/2(1)$ & 0.524 & 1499 & 0.66\\\hline\end{tabular}
\noindent\begin{tabular}{|l||r|r|r|}\hline& $a$ & $b$ & $c$\\\hline\hline$47/2(2)$ & 0.509 & 642 & 1.57\\\hline$47/2(3)$ & 0.562 & 2105 & 0.45\\\hline$49/2(1)$ & 0.542 & 1023 & 1.01\\\hline$49/2(2)$ & 0.390 & 3004 & 0.36\\\hline$49/2(3)$ & 0.561 & 1841 & 0.53\\\hline$49/2(4)$ & 0.453 & 400 & 2.14\\\hline$51/2(1)$ & 0.606 & 807 & 1.39\\\hline$51/2(2)$ & 0.537 & 2644 & 0.36 \\\hline$51/2(3)$ & 0.519 & 889 & 1.10\\\hline$53/2(1)$ & 0.486 & 2725 & 0.37\\\hline$53/2(2)$ & 0.426 & 2126 & 0.48\\\hline$53/2(3)$ & 0.630 & 986 & 1.13\\\hline$53/2(4)$ & 0.417 & 1304 & 0.72\\\hline$55/2(1)$ & 0.568 & 1165 & 0.88 \\\hline$55/2(2)$ & 0.516 & 1287 & 0.77\\\hline$55/2(3)$ & 0.670 & 216 & 7.00 \\\hline$55/2(4)$ & 0.463 & 3322 & 0.31\\\hline$57/2(1)$ & 0.568 & 1844 & 0.53\\\hline$57/2(2)$ & 0.529 & 438 & 2.36\\\hline$57/2(3)$ & 0.438 & 2903 & 0.36\\\hline$57/2(4)$ & 0.530 & 943 & 1.08\\\hline$59/2(1)$ & 0.568 & 1060 & 0.98\\\hline$59/2(2)$ & 0.496 & 1539 & 0.64\\\hline$59/2(3)$ & 0.586 & 1068 & 0.99\\\hline$59/2(4)$ & 0.416 & 3728 & 0.30\\\hline$61/2(1)$ & 0.562 & 1818 & 0.54\\\hline$61/2(2)$ & 0.430 & 3023 & 0.35\\\hline$61/2(3)$ & 0.579 & 199 & 6.74\\\hline$61/2(4)$ & 0.502 & 1848 & 0.54\\\hline$61/2(5)$ & 0.406 & 472 & 1.61\\\hline\end{tabular}\end{multicols}

\newpage\noindent Best fit parameters (rounded) for the Laplace distribution for all examples with $10^7$ coefficients:\begin{multicols}{3}\noindent\begin{tabular}{|l||r|r|}\hline& $b$ & $c$ \\\hline\hline$13/2$ & 1172 & 0.908 \\\hline$17/2$ & 1503 & 0.683\\\hline$19/2$ & 2600 & 0.388\\\hline$21/2$ & 1740 & 0.592\\\hline$23/2$ & 1220 & 0.850\\\hline$25/2(1)$ & 1260 & 0.791\\\hline$25/2(2)$ & 1884 & 0.511\\\hline$27/2$ & 1923 & 0.542\\\hline$29/2$ & 725 & 1.317\\\hline$31/2(1)$ & 1914 & 0.533\\\hline$31/2(2)$ & 1258 & 0.789\\\hline$33/2(1)$ & 1302 & 0.793\\\hline$33/2(2)$ & 2202 & 0.444\\\hline$35/2(1)$ & 3147 & 0.295\\\hline$35/2(2)$ & 62 & 17.400 \\\hline$37/2(1)$ & 2140 & 0.440\\\hline$37/2(2)$ & 1930 & 0.522\\\hline$37/2(3)$ & 429 & 2.341\\\hline$39/2(1)$ & 1495 & 0.669\\\hline$39/2(2)$ & 1881 & 0.550\\\hline\end{tabular}
\noindent\begin{tabular}{|l||r|r|}\hline& $b$ & $c$ \\\hline\hline$41/2(1)$ & 1367 & 0.740\\\hline$41/2(2)$ & 2458 & 0.395\\\hline$41/2(3)$ & 1163 & 0.897\\\hline$43/2(1)$ & 981 & 0.980\\\hline$43/2(2)$ & 1160 & 0.899\\\hline$43/2(3)$ & 3021 & 0.309\\\hline$45/2(1)$ & 2461 & 0.409\\\hline$45/2(2)$ & 774 & 1.331\\\hline$45/2(3)$ & 1366 & 0.699\\\hline$47/2(1)$ & 1532 & 0.655\\\hline$47/2(2)$ & 647 & 1.539\\\hline$47/2(3)$ & 2219 & 0.460\\\hline$49/2(1)$ & 1063 & 0.953\\\hline$49/2(2)$ & 2650 & 0.333\\\hline$49/2(3)$ & 1940 & 0.527\\\hline$49/2(4)$ & 381 & 2.493\\\hline$51/2(1)$ & 871 & 1.196 \\\hline$51/2(2)$ & 2734 & 0.370 \\\hline$51/2(3)$ & 903 & 1.077 \\\hline$53/2(1)$ & 2687 & 0.364\\\hline\end{tabular}
\noindent\begin{tabular}{|l||r|r|}\hline& $b$ & $c$ \\\hline\hline$53/2(2)$ & 1958 & 0.479\\\hline$53/2(3)$ & 1090 & 0.962\\\hline$53/2(4)$ & 1188 & 0.778\\\hline$55/2(1)$ & 1233 & 0.830 \\\hline$55/2(2)$ & 1306 & 0.760 \\\hline$55/2(3)$ & 225 & 4.846 \\\hline$55/2(4)$ & 3197 & 0.303 \\\hline$57/2(1)$ & 1954 & 0.525\\\hline$57/2(2)$ & 448 & 2.206\\\hline$57/2(3)$ & 2716 & 0.349\\\hline$57/2(4)$ & 969 & 1.039\\\hline$59/2(1)$ & 1121 & 0.912 \\\hline$59/2(2)$ & 1533 & 0.642 \\\hline$59/2(3)$ & 1144 & 0.902 \\\hline$59/2(4)$ & 3389 & 0.274 \\\hline$61/2(1)$ & 1918 & 0.532\\\hline$61/2(2)$ & 2796 & 0.336\\\hline$61/2(3)$ & 211 & 4.830\\\hline$61/2(4)$ & 1852 & 0.536\\\hline$61/2(5)$ & 423 & 2.131\\\hline\end{tabular}\end{multicols}
\newpage\noindent Best fit parameters (rounded) for the Cauchy distribution for all examples with $10^7$ coefficients:\begin{multicols}{2}\noindent \begin{tabular}{|l||r|r|r|}\hline& $a$ & $b$ & $c$\\\hline\hline$13/2$ & 143 & 0.14 & 0.49\\\hline$17/2$ & 154 & 0.12 & 0.61\\\hline$19/2$ & 152 & 0.07 & 0.82\\\hline$21/2$ & 150 & 0.10 & 0.64\\\hline$23/2$ & 156 & 0.15 & 0.54\\\hline$25/2(1)$ & 163 & 0.15 & 0.59\\\hline$25/2(2)$ & 169 & 0.10 & 0.78\\\hline$27/2$ & 148 & 0.09 & 0.66\\\hline$29/2$ & 162 & 0.25 & 0.46\\\hline$31/2(1)$ & 156 & 0.09 & 0.68\\\hline$31/2(2)$ & 163 & 0.15 & 0.59\\\hline$33/2(1)$ & 156 & 0.14 & 0.56\\\hline$33/2(2)$ & 163 & 0.08 & 0.81\\\hline$35/2(1)$ & 212 & 0.07 & 1.17\\\hline$35/2(2)$ & 422 & 7.75 & 0.19 \\\hline$37/2(1)$ & 200 & 0.10 & 0.93\\\hline$37/2(2)$ & 156 & 0.09 & 0.71\\\hline$37/2(3)$ & 5310 & 13.99 & 1.96\\\hline$39/2(1)$ & 175 & 0.13 & 0.67\\\hline$39/2(2)$ & 145 & 0.09 & 0.66\\\hline$41/2(1)$ & 161 & 0.13 & 0.61\\\hline$41/2(2)$ & 196 & 0.09 & 0.95\\\hline$41/2(3)$ & 147 & 0.14 & 0.53\\\hline$43/2(1)$ & 192 & 0.22 & 0.60\\\hline$43/2(2)$ & 145 & 0.14 & 0.51\\\hline$43/2(3)$ & 243 & 0.09 & 1.23\\\hline$45/2(1)$ & 287 & 0.13 & 1.10\\\hline$45/2(2)$ & 153 & 0.24 & 0.44\\\hline$45/2(3)$ & 206 & 0.17 & 0.74\\\hline$47/2(1)$ & 164 & 0.12 & 0.66\\\hline\end{tabular}
\noindent \begin{tabular}{|l||r|r|r|}\hline& $a$ & $b$ & $c$\\\hline\hline$47/2(2)$ & 165 & 0.29 & 0.43\\\hline$47/2(3)$ & 151 & 0.08 & 0.74\\\hline$49/2(1)$ & 163 & 0.17 & 0.54\\\hline$49/2(2)$ & 211 & 0.09 & 1.14\\\hline$49/2(3)$ & 152 & 0.09 & 0.70\\\hline$49/2(4)$ & 1070 & 3.14 & -0.90\\\hline$51/2(1)$ & 13 & 0.02 & 0.13\\\hline$51/2(2)$ & 12 & 0.00 & 0.23\\\hline$51/2(3)$ & 952 & 1.20 & 1.23 \\\hline$53/2(1)$ & 169 & 0.07 & 0.91\\\hline$53/2(2)$ & 219 & 0.12 & 0.94\\\hline$53/2(3)$ & 144 & 0.15 & 0.49\\\hline$53/2(4)$ & 216 & 0.20 & 0.73\\\hline$55/2(1)$ & 580 & 0.54 & 1.06 \\\hline$55/2(2)$ & 1857 & 1.61 & 2.04 \\\hline$55/2(3)$ & 1459 & 7.01 & 0.82 \\\hline$55/2(4)$ & 3681 & 1.29 & 4.70 \\\hline$57/2(1)$ & 149 & 0.59 & 0.69\\\hline$57/2(2)$ & 65049 & 164.83 & 7.08\\\hline$57/2(3)$ & 208 & 0.09 & 1.05\\\hline$57/2(4)$ & 163 & 0.19 & 0.52\\\hline$59/2(1)$ & 1383 & 1.41 & 1.57 \\\hline$59/2(2)$ & 2525 & 1.86 & 2.62 \\\hline$59/2(3)$ & 385 & 28.90 & 24.97 \\\hline$59/2(4)$ & 88133 & 0.30 & 2.55 \\\hline$61/2(1)$ & 151 & 0.09 & 0.69\\\hline$61/2(2)$ & 242 & 0.10 & 1.17\\\hline$61/2(3)$ & 3403 & 18.40 & -1.07\\\hline$61/2(4)$ & 169 & 0.10 & 0.74\\\hline$61/2(5)$ & 1346 & 3.52 & -1.12\\\hline\end{tabular}\end{multicols}
\newpage\noindent RMS values (rounded) for all examples with $10^7$ coefficients:\begin{multicols}{2}\noindent\begin{tabular}{|l||r|r|r|r|}\hline & GG & GGG & Laplace & Cauchy \\\hline\hline$13/2$ & 19 & 18 & 39 & 33 \\\hline$17/2$ & 21 & 18 & 29 & 28 \\\hline$19/2$ & 27 & 19 & 30 & 34 \\\hline$21/2$ & 21 & 18 & 31 & 33 \\\hline$23/2$ & 19 & 18 & 28 & 26 \\\hline$25/2(1)$ & 21 & 18 & 22 & 22 \\\hline$25/2(2)$ & 33 & 19 & 35 & 25 \\\hline$27/2$ & 19 & 17 & 37 & 41 \\\hline$29/2$ & 18 & 16 & 18 & 20\\\hline$31/2(1)$ & 22 & 19 & 28 & 33 \\\hline$31/2(2)$ & 21 & 17 & 22 & 21 \\\hline$33/2(1)$ & 19 & 18 & 28 & 27 \\\hline$33/2(2)$ & 27 & 18 & 28 & 30 \\\hline$35/2(1)$ & 48 & 20 & 57 & 40 \\\hline$35/2(2)$ & 5 & 5 & 5 & 5 \\\hline$37/2(1)$ & 26 & 14 & 30 & 25 \\\hline$37/2(2)$ & 20 & 15 & 22 & 25 \\\hline$37/2(3)$ & 9 & 9 & 9 & 10 \\\hline$39/2(1)$ & 16 & 12 & 16 & 20 \\\hline$39/2(2)$ & 17 & 15 & 28 & 34 \\\hline$41/2(1)$ & 14 & 12 & 16 & 20 \\\hline$41/2(2)$ & 22 & 14 & 24 & 27 \\\hline$41/2(3)$ & 13 & 11 & 14 & 18 \\\hline$43/2(1)$ & 13 & 9 & 13 & 14 \\\hline$43/2(2)$ & 14 & 13 & 22 & 25 \\\hline$43/2(3)$ & 41 & 19 & 50 & 41 \\\hline$45/2(1)$ & 21 & 16 & 23 & 31 \\\hline$45/2(2)$ & 11 & 10 & 13 & 16 \\\hline$45/2(3)$ & 15 & 11 & 17 & 17 \\\hline$47/2(1)$ & 16 & 13 & 16 & 21 \\\hline\end{tabular}
\noindent 
\begin{tabular}{|l||r|r|r|r|}
\hline 
& GG & GGG & Laplace & Cauchy \\\hline\hline$47/2(2)$ & 10 & 9 & 10 & 12 \\\hline$47/2(3)$ & 19 & 16 & 25 & 34 \\\hline$49/2(1)$ & 13 & 11 & 14 & 17 \\\hline$49/2(2)$ & 59 & 24 & 72 & 49 \\\hline$49/2(3)$ & 17 & 15 & 23 & 31 \\\hline$49/2(4)$ & 10 & 10 & 10 & 10 \\\hline$51/2(1)$ & 18 & 17 & 24 & 20 \\\hline$51/2(2)$ & 23 & 18 & 27 & 38 \\\hline$51/2(3)$ & 20 & 17 & 20 & 21 \\\hline$53/2(1)$ & 25 & 16 & 25 & 28 \\\hline$53/2(2)$ & 20 & 12 & 26 & 25 \\\hline$53/2(3)$ & 13 & 13 & 23 & 26 \\\hline$53/2(4)$ & 17 & 10 & 20 & 17 \\\hline$55/2(1)$ & 19 & 17 & 24 & 23 \\\hline$55/2(2)$ & 21 & 18 & 22 & 23 \\\hline$55/2(3)$ & 12 & 12 & 13 & 13 \\\hline$55/2(4)$ & 30 & 18 & 34 & 39 \\\hline$57/2(1)$ & 18 & 16 & 24 & 32\\\hline$57/2(2)$ & 12 & 12 & 12 & 12 \\\hline$57/2(3)$ & 27 & 15 & 33 & 31 \\\hline$57/2(4)$ & 13 & 11 & 14 & 16 \\\hline$59/2(1)$ & 19 & 18 & 24 & 22 \\\hline$59/2(2)$ & 23 & 17 & 23 & 25 \\\hline$59/2(3)$ & 19 & 18 & 26 & 24 \\\hline$59/2(4)$ & 40 & 19 & 54 & 48 \\\hline$61/2(1)$ & 17 & 15 & 23 & 29 \\\hline$61/2(2)$ & 34 & 18 & 43 & 42 \\\hline$61/2(3)$ & 9 & 9 & 9 & 9 \\\hline$61/2(4)$ & 17 & 13 & 17 & 23 \\\hline$61/2(5)$ & 11 & 9 & 13 & 12 \\\hline
\end{tabular}
\end{multicols}

\end{document}